\documentclass[3p, authoryear, referee]{elsarticle}
\usepackage{url}
\usepackage{hyperref}
\hypersetup{
    colorlinks=true,
    linkcolor=blue,
    filecolor=blue,      
    urlcolor=blue,
}

\usepackage{natbib}
\usepackage{makecell}
\setcellgapes{7pt}

\usepackage{graphics,graphicx,multirow,dsfont}
\usepackage{amsmath,amssymb,verbatim,epsfig, amsthm}
\usepackage{booktabs}
\usepackage{multirow}
\usepackage{color, colortbl}
\usepackage{bbold}
\usepackage{subfigure}
\usepackage[ruled,vlined,plain,noend]{algorithm2e}
\usepackage{ulem}

\newtheorem{theorem}{Theorem}
\newtheorem{proposition}{Proposition}

\newtheorem{corollary}{Corollary}
\newtheorem{defin}{\bf Definition}

\def\R{\mathds{R}}
\def\X{\mathds{X}}
\def\e{\mathrm{e}}

\def\E{\mathds{E}}

\def\d{\mathrm{d}}

\def\bm{{\bf m}}

\def\simiid{\stackrel{\mbox{iid}}{\sim}}

\allowdisplaybreaks

\usepackage{times}
\usepackage{bm,mathrsfs,dsfont,color}
\usepackage{graphicx}
\allowdisplaybreaks

\usepackage{multirow}

\newcommand{\indic}{\mathds{1}}

\newcommand{\ddr}{\mathrm{d}}
\newcommand{\edr}{\mathrm{e}}

\def\simiid{\stackrel{\mbox{\scriptsize{iid}}}{\sim}}

\makeatletter
\def\ps@pprintTitle{%
 \let\@oddhead\@empty
 \let\@evenhead\@empty
 \def\@oddfoot{}%
 \let\@evenfoot\@oddfoot}
\makeatother
 
\begin{document}

\title{
	Inner spike and slab Bayesian nonparametric models}

 \author{Antonio Canale}
 \address{Department of Statistical Sciences, 
University of Padova\\
Via Cesare Battisti, 241/243
35121 Padova, Italy}
 
  \author{Antonio Lijoi}
 \address{Department of Decision Sciences and BIDSA, Bocconi University\\
 via R\"ontgen 1, 20136 Milan, Italy }
 
 \author{Bernardo Nipoti\corref{cor1}}
 \address{Department of Economics, Management and Statistics, University of Milano-Bicocca\\ Piazza dell'Ateneo Nuovo, 1, 20126, Milan, Italy}
  \ead{bernardo.nipoti@unimib.it}
 
  \author{Igor Pr\"unster}
 \address{Department of Decision Sciences and BIDSA, Bocconi University\\[3pt]
 via R\"ontgen 1, 20136 Milan, Italy }
 
  \cortext[cor1]{Corresponding author}

\begin{abstract}
Discrete Bayesian nonparametric models whose expectation is a convex linear combination of a point mass at some point of the support and a diffuse probability distribution 
allow to incorporate strong prior information, while still being extremely flexible.
Recent contributions in the statistical literature have successfully implemented such a modelling strategy in a variety of applications, including density estimation, nonparametric regression and model-based clustering. We provide a thorough study of a large class of nonparametric models we call \textit{inner spike and slab hNRMI models}, which are obtained by considering homogeneous normalized random measures with independent increments (hNRMI) with base measure given by a convex linear combination of a point mass and a diffuse probability distribution. In this paper we investigate the distributional properties of these models and our results include: i) the exchangeable partition probability function they induce, ii) the distribution of the number of distinct values in an exchangeable sample, iii) the posterior predictive distribution, and iv) the distribution of the number of elements that coincide with the only point of the support with positive probability. Our findings are the main building block for an actual implementation of Bayesian inner spike and slab hNRMI models by means of a generalized P\'olya urn scheme.
\end{abstract}


\maketitle
\section{Introduction}

One of the most appealing aspects of Bayesian nonparametric modelling is its ability to flexibly account for a rich  variety of patterns beyond those described by specific parametric models. In many applications, however, the need of flexibility comes along with the availability of valuable prior information on some attributes of the data. For example, when monitoring a process which varies over time, interest typically lies in
discriminating between a known baseline behaviour and departures from the same. Similarly, when modelling functional data, it is often the case that some of the observations are expected to exhibit specific regular features while others might display more erratic traits. Both these situations call for a model which might be able to accurately capture the expected baseline behaviour while, at the same time, being flexible in dealing with the more irregular observations. Many recent contributions have formalized this idea by considering nonparametric mixture models defined so that the expected baseline behaviour is assigned positive prior probability. We review here two simple strategies that allow to achieve this goal, henceforth referred to as  \textit{inner} and \textit{outer} spike and slab models, following a terminology introduced in \citet{Can17}. 

We consider an observation $X$ taking values in some space $\X$ and we suppose that the baseline behaviour is described by a point $x_0  \in\X$. The outer spike and slab model is a convex linear combination of a Dirac delta at $x_0$ and a discrete nonparametric measure on $\X$. More specifically, the model
assumes that $X \sim \tilde Q$, where $\tilde Q$ on $\X$ is defined as 
\begin{equation}
{\tilde Q}= \zeta \delta_{x_0} + (1-\zeta) Q^*,
\label{eq:outer}
\end{equation}
where $\zeta \in (0,1)$ represents the probability of the event $\{X=x_0\}$, and $Q^*$ is a discrete random probability measure such that $\E[Q^*]=P^*$, with $P^*$ a diffuse probability measure. Consequently, one has $\E[\tilde Q]= \zeta \delta_{x_0} + (1-\zeta) P^*$. 
Alternatively, the inner spike and slab model is a discrete random probability measure $\tilde P=\sum_{j}\tilde p_j\,\delta_{Z_j}$ with random locations $Z_j$'s on $\X$ specified in such a way that 
\begin{equation}
\E[\tilde P]=\zeta \delta_{x_0}+(1-\zeta)P^*,
\label{eq:inner}
\end{equation}
where $\zeta \in (0,1)$ can be interpreted, again, as the probability that $X$ coincides with $x_0$ and $P^*$ is a diffuse probability distribution on $\X$. 
The two models are characterized by analogous spike and slab structures \citep{mitchell(88)} and are such that $\E[\tilde Q]=\E[\tilde P]=\zeta \delta_{x_0}+(1-\zeta)P^*$. Nonetheless,  the outer and inner specifications differ in that the spike at $x_0$ is external to the nonparametric component of the model, in the former, as displayed by \eqref{eq:outer}, while it is included in the nonparametric part, in the latter, with $x_0$ being the only atomic component of the otherwise nonatomic expected value \eqref{eq:inner}. Inner and outer spike and slab specifications have been largely adopted in the recent literature. 
	\citet{Cas19} model pneumonia and influenza mortality in space and time by adopting the outer spike and slab mixture approach described in \eqref{eq:outer}, with $Q^*$ distributed as a Dirichlet process (DP) \citep{Ferg73}. Similarly, \citet{Dangelo21} use the inner spike and slab approach to Bayesian nonparametric mixtures 
	to analyze data on neural activity of animals, measured via calcium imaging. Another example is provided by 
	\citet{scar:duns:2009} and \citet{Can17}, where the same functional data set on women fertility, is analysed, respectively, by means of outer and inner spike and slab mixture models defined on functional spaces.
Similar approaches have been used also in the context of variable selection. \citet{Dun(12)}, \citet{Mac(07)}, \citet{Yan(12)} and \cite{deiorio(16)}, for example,  adopt a DP specified such that \eqref{eq:inner} holds with $x_0 = 0$, to simultaneously allow for variable selection and clustering of the variables. The same construction is adopted in \citet{Sua(16)} to model wavelet coefficients of functional data so  to induce sparsity.  Models to infer differential gene expression based on a DP with inner spike and slab base measure, can be found in \cite{peter05} and \citet{	guindani2014bayesian}, while applications to multiple testing problems are proposed by \cite{surgio} and \cite{dahl09}.   

Most of the aforementioned contributions consider inner spike and slab models based on a DP.  In fact, when assuming $\tilde P$ is distributed as a DP, specifications leading to a spike at $x_0$ in \eqref{eq:inner} 
do not change the structure of the resulting predictive distribution thanks to the conjugacy of the DP \citep{Ferg73}. However, when $\tilde P$ is not a DP, specifications implying 
an expected value with atoms, such as the spike at $x_0$ in  \eqref{eq:inner}, considerably change the posterior predictive structure of the process \citep[see][]{San(06),Can17}, thus implying challenging technical issues that need to be addressed in order to perform Bayesian posterior inference. On the other hand, working with outer spike and slab models is in general less cumbersome as the nonparametric component of the model has diffuse base measure, for which, thus, standard techniques can be used. 

The inner spike and slab model has been thoroughly studied by \citet{Can17} in the case of $\tilde P$ being distributed as a Pitman--Yor (PY) process \citep{Per(92),PY1997}. More recently, \cite{Bas20} investigated the asymptotic behaviour of species sampling models with non-diffuse base measure. In this paper we focus on the flexible class of nonparametric priors obtained by normalizing homogeneous 
completely random measures (CRMs) \citep{Reg03}, henceforth called 
homogeneous 
normalized random measures with independent increments (hNRMIs). Their stick-breaking representation has been derived in \citet{Fav16}.  We study the distributional properties of hNRMIs with spike and slab specification \eqref{eq:inner}, thus providing the essential building block for carrying out posterior inference with inner spike and slab hNRMI models. Our investigations also highlight that, besides the apparent similarities, as far as hNRMIs are concerned, the inner and the outer spike and slab models are structurally different. While having coinciding prior expectations, the two model specifications differ in terms of prior variability, with the variance turning out to be larger for the inner spike and slab model. In this sense, the inner spike and slab prior is less informative than the outer one. Our investigation on this point is underpinned by a numerical study, whose results are in line with the findings of \citet{Can17} for the PY case, where the inner approach is showed to appear more robust than the outer one, when dealing with misspecified priors.

The paper is organised as follows. In Section \ref{sec:nrmi} we provide a succinct introduction to hNRMIs and compare the prior variance of functionals of inner and outer spike and slab models. Section \ref{sec:main} presents a thorough investigation of the distributional properties of the inner spike and slab hNRMI model. As notable examples, in Section \ref{sec:stable}, the main results of the paper are displayed for the special cases of $\sigma$-stable and normalized inverse Gaussian hNRMIs. Section \ref{sec:algo} provides a generalised P\'olya urn scheme for the inner spike and slab hNRMI model. This is then implemented, for the $\sigma$-stable case, to allow for a numerical comparison of inner and outer spike and slab models, both \textit{a priori} and \textit{a posteriori}. Finally, the proofs of the main results are reported in Section \ref{sec:proofs}.

\section{Homogeneous normalized random measures and spike and slab models} \label{sec:nrmi}
We concisely recall the basics of CRMs and hNRMIs tailored to the present contribution. For an account on their role in Bayesian nonparametric statistics one can refer to \citet{Lp10}. 
	A random measure $\tilde\mu$ on $\X$ is said \textit{completely random} if for any collection of pairwise disjoint measurable subsets $A_1,\ldots,A_k$ of $\X$, and for any $k\ge 2$, the random variables $\tilde\mu(A_1),\ldots,\tilde\mu(A_k)$ are mutually independent.  If $\tilde\mu$ is without fixed points of discontinuity, which is assumed throughout the paper, for any function $f:\X\to\R^+$ one has 
	\begin{equation*}
	\E\left[\edr^{-\int_\X f(x)\tilde\mu(\ddr x)}\right]=
	\exp\left\{-\int_{\R^+\times\X}\left(1-\edr^{-s f(x)}\right)\nu(\ddr s,\ddr x) \right\}
	\label{eq:laplace_expon}
	\end{equation*}
	where $\nu$, a measure on $\R^+\times \X$, is the L\'evy intensity characterizing $\tilde \mu$.  For the ease of illustration, 
	henceforth we consider the class of almost surely finite homogeneous CRMs, corresponding to L\'evy intensities that admit a factorization of the type 
	$
	\nu(\ddr s,\ddr x)=\rho(s)\ddr s cP_0(\ddr x)
	$
	for some measurable function $\rho:\R^+\to\R^+$ and some constant $c>0$. Noteworthy examples that will be considered in this paper are the $\sigma$-stable process and the inverse Gaussian process, characterized respectively by $\rho(s)=\sigma s^{-1-\sigma} /\Gamma(1-\sigma)$, for some $\sigma\in(0,1)$, and $\rho(s)=s^{-3/2}\edr^{-\tau s}/(2\sqrt{\pi})$, for some $\tau>0$. A random probability measure is then obtained by normalization as $\tilde P=\tilde\mu/\tilde\mu(\X)$ and denoted by $\tilde P\sim\mbox{hNRMI}(\rho;c,P_0)$. The study of this class of nonparametric priors was first considered in \cite{Reg03}. The special cases considered in this paper are two popular nonparametric priors, namely the $\sigma$-stable hNRMI and the normalized inverse Gaussian (N-IG) process \citep{Lij05}, which are obtained by normalizing a $\sigma$-stable CRM and an inverse Gaussian process, respectively. See e.g. \cite{Lp10} and \cite{statsci13} for a review of their inferential properties.
	
	We are now in the position to highlight a structural difference between inner and outer hNRMI spike and slab models. Specifically, we let $\tilde P\sim \text{hNRMI}(\rho;c,P_0)$ with spike and slab base measure 
		\begin{equation}\label{eq:base_inner}
			P_0=\zeta \delta_{x_0}+(1-\zeta)P^*.
		\end{equation} 
	This defines an inner spike and slab model and implies $\E[\tilde P]=P_0$ \citep{Jam06}. 	Moreover, we consider an outer spike and slab model $\tilde Q$, defined as in \eqref{eq:outer} with $Q^*\sim \text{hNRMI}(\rho;c,P^*)$. While the two random probability measures have the same expectation, namely $\zeta \delta_{x_0}+(1-\zeta)P^*$, they differ in terms of variance, as displayed in the next proposition. Before stating the result, we define $\tau_q(u)=\int_0^\infty s^q \edr^{-us} \rho(s) \ddr s$, for any integer $q\ge 1$, and $\psi(u)=\int_0^\infty(1-\edr^{-u s})\rho(s) \ddr s$, for any $u>0$. Moreover, for any probability measure $P$, we use the notation $P(f)=\int f\,\ddr P$.\medskip
	
	\begin{proposition}\label{prop:var}
	Let $f:\X\to\R$ any measurable function such that $P^*(f^2)<\infty$. Then 
\begin{equation*}\label{eq:var}
    \text{var}(\tilde P(f) )- \text{var}(\tilde Q(f))=p\,\zeta(1-\zeta)\,P^*\Big( \{f-f(x_0)\}^2\Big) \geq 0,
\end{equation*}
where $p=c\int_0^\infty u\edr^{-c\psi(u)}\tau_{2}(u)\ddr u\in(0,1)$.
	\end{proposition}
	Proposition \ref{prop:var} indicates that, for the class of hNRMIs, the inner spike and slab model is characterized by larger prior uncertainty than the outer spike and slab model with the same expectation. In other terms, while both models are centered at the same prior guess, the inner specification of the model is less informative. This aspect will be further investigated in the numerical study presented in Section \ref{sec:algo}.
	
	\section{Main results}\label{sec:main}
	
We now investigate the distributional properties of $\tilde P\sim\text{hNRMI}(\rho;c,P_0)$, where $P_0$ is a spike and slab base measure defined as in \eqref{eq:base_inner}. Our results include: i) the exchangeable partition probability function (EPPF) induced by $\tilde P$, that is the probability of observing a specific sample displaying $k$ distinct values, henceforth also referred to as clusters, with corresponding frequencies summarized by the vector $(n_1,\ldots,n_k)$; ii) the distribution of $K_n$, the number of distinct values in an exchangeable sample $X^{(n)}=(X_1,\ldots,X_n)$ such that $X_i\mid \tilde P \simiid \tilde P$; iii) the predictive distribution for one observation in the exchangeable sample, conditionally on the observation of the others; iv) the distribution of $N^{(n)}_0$, that is the number of elements of $X^{(n)}$ that coincide with $x_0$. Henceforth we use the notation $\Pi_{k}^{(n)}(n_1,\ldots,n_k;\zeta)$ to denote the EPPF induced by a hNRMI with spike and slab base measure \eqref{eq:base_inner} parametrized by $\zeta\in(0,1)$. Similarly, we will write $\Pr(\,\cdot\,;\zeta)$ to stress the fact that the probability of an event depends on the parameter $\zeta$. Accordingly, the notation $\Pi_{k}^{(n)}(n_1,\ldots,n_k;0)$ and $\Pr(\,\cdot\,;0)$ will refer to the EPPF and the probability measure induced by a hNRMI with diffuse base measure $P^*$. 
Finally, for any positive integer $n$, and any  $i=1,\ldots,n$, we define the function
\begin{equation}\label{eq:xi}
\xi_{n,i}(u)=\frac{1}{i!}\sum_{j=0}^i (-1)^{n-j} \binom{i}{j}\psi^{i-j}(u)\frac{\ddr^n}{\ddr u^n}\left[\psi^j(u)\right].
\end{equation}

\begin{theorem}\label{eppf:nrmi}
If $X_i\mid \tilde P\simiid \tilde P$, for $i=1,\ldots,n$ and for any $n\ge 1$, and $\tilde P\sim\mbox{\rm hNRMI}(\rho;c,P_0)$, with $P_0$ as in \eqref{eq:base_inner}, the EPPF induced by $\tilde P$ is 
\begin{multline}\label{eq:eppf:nrmi}
	\Pi_k^{(n)}(n_1,n_2,\ldots,n_k;\zeta)=\frac{1}{\Gamma(n)}\left\{c^k(1-\zeta)^k\int_0^\infty u^{n-1}\edr^{-c\psi(u)}\prod_{m=1}^{k}\tau_{n_m}(u)\ddr u\right.\\
	\left.+c^{k-1}(1-\zeta)^{k-1}\sum_{l=1}^{k}\sum_{i=1}^{n_l}c^{i}\zeta^{i}\int_{0}^\infty u^{n-1}\edr^{-c\psi(u)}\xi_{n_{l,i}}(u)\prod_{m\neq l}\tau_{n_m}(u)\ddr u\right\},
	\end{multline}
where we agree that $\prod_{m\ne \ell}\tau_{n_m}(u)\equiv 1$ when $k=1$.
\end{theorem}

Out of the $k+1$ summands on the right-hand side of \eqref{eq:eppf:nrmi}, the first one refers to the case in which none of the blocks of the partition coincide with $x_0$, while the remaining $k$ terms account for the cases in which the $l$-th cluster is identified by the atom $x_0$, for $l=1,2,\ldots,k$.\\
We recall that the EPPF of a hNRMI with diffuse base measure $P^*$ equals
\begin{equation}
\label{eq:eppf_hnrmi}
\Pi_k^{(n)}(n_1,\ldots,n_k;0)=\frac{c^k}{\Gamma(n)}\int_0^\infty u^{n-1}\edr^{-c\psi(u)}\prod_{j=1}^k\tau_{n_j}(u)\ddr u.
\end{equation}
Henceforth, for any $n_j\in\{1,\ldots,n-k+1\}$, let 
$\mathcal{X}_{k-1,0}^{(n-n_j)}(\bm{n}_{-j})$ denote a sample of size $n-n_j$ clustered into $k-1$ groups with respective frequencies $\bm{n}_{-j}=(n_1,\ldots,n_{j-1},n_{j+1},\ldots,n_k)$ such that: (i) none of its elements equals $x_0$ and (ii) when extended to a sample of size $n$, the overall number of observations that coincide with $x_0$ equals $n_j$, namely $N_{0}^{(n)}=\mbox{card}\{i:\: X_i=x_0\}=n_j$.
One can then state the following result.
\begin{corollary}\label{thm:interpret_eppf_nrmi}
	If $X_i\mid \tilde P\simiid \tilde P$, for $i=1,\ldots,n$ and for any $n\ge 1$, and $\tilde P\sim\mbox{\rm hNRMI}(\rho;c,P_0)$, with $P_0$ as in \eqref{eq:base_inner}, the EPPF induced by $\tilde P$ is
	\begin{multline}\label{eq:eppf_bis}
	\Pi_k^{(n)}(n_1,\ldots,n_k;\zeta)=
	(1-\zeta)^k\Pi_k^{(n)}(n_1,\ldots,n_k;0)\\[4pt]
	+(1-\zeta)^{k-1}\sum_{j=1}^k \Pi_{k-1}^{(n-n_j)}(n_1,\ldots,n_{j-1},n_{j+1},\ldots,n_k;0)
	\Pr (N_0^{(n)}=n_j\mid \mathcal{X}_{k-1,0}^{(n-n_j)}(\bm{n}_{-j}) 
	;\zeta ).
	\end{multline}
\end{corollary}

It is apparent that the first summand in \eqref{eq:eppf_bis} accounts for the case where none of the $k$ clusters is identified by $x_0$, so one has the standard EPPF that corresponds to a diffuse base measure. On the other hand, the second summand takes into account the possibility that $x_0$ identifies one of the $k$ clusters in the partition. 

Starting from Theorem \ref{thm:interpret_eppf_nrmi} we obtain the distribution of $K_n$, the number of distinct values in the sample $X^{(n)}$, as reported in the next result. To this end, for any $a,u>0$ we will denote with $X_{a,u}$  a random variable such that its density function is $f_{a,u}(s)\propto s^a\,\edr^{-us}\,\rho(s)\,\indic_{(0,+\infty)}(s)$. It can be seen that if the following condition holds true
\begin{itemize}
	\item[(H1)] $\rho$ is such that for any finite collection of independent random variables $X_{a_1,u},\ldots,X_{a_n,u}$, the distribution of $\sum_{i=1}^n X_{a_i,u}$ depends on $(a_1,\ldots,a_n)$ only through $\sum_{i=1}^n a_i$, 
\end{itemize}
then $\mbox{Pr}(N_0^{(n)}=r\mid\mathcal{X}_{k-1,0}^{(n-r)}(\bm{n}_{-j});\zeta)=\mbox{Pr}(N_0^{(n)}=r\mid K_{n-r}=k-1;\zeta)$. This is used to prove the following theorem. \medskip  

\begin{theorem}\label{thm:K_n}
Let $X_i\mid \tilde P\simiid \tilde P$, for $i=1,\ldots,n$ and for any $n\ge 1$, and $\tilde P\sim\mbox{\rm hNRMI}(\rho;c,P_0)$, with $P_0$ as in \eqref{eq:base_inner}. Moreover, $\rho$ is such that {\rm (H1)} holds true. Then, for any $k\in\{1,\ldots,n\}$, the distribution of $K_n$, 
parametrized by $\zeta\in(0,1)$,  
is given by
\begin{multline*}\label{eq:Kn}
\Pr(K_n=k;\zeta)=(1-\zeta)^k \Pr(K_n=k;0)\\
+(1-\zeta)^{k-1}\sum_{r=1}^{n-k+1}\binom{n}{r}  \Pr(K_{n-r}=k-1;0) 
	\Pr(N_0^{(n)}=r\mid K_{n-r}=k-1 
	;\zeta ).
\end{multline*}
\end{theorem}

It is worth stressing that for all the examples we consider henceforth condition (H1) holds true. In view of these findings, the predictive distributions associated to the sequence $(X_i)_{i\ge 1}$ can now be easily determined.  
We suppose that the observed sample $X^{(n)}$ displays $k$ distinct values $x_1^*,\ldots,x_k^*$, with respective frequencies $n_1,\ldots,n_k$, and state the following theorem. \medskip

\begin{theorem}\label{thm:predictive}
If $X_i\mid \tilde P\simiid \tilde P$, for $i=1,\ldots,n$ and for any $n\ge 1$, and $\tilde P\sim\mbox{\rm hNRMI}(\rho;c,P_0)$, with $P_0$ as in \eqref{eq:base_inner}, the predictive distribution of $X_{n+1}$ conditionally on $X^{(n)}$ is
\begin{itemize}
    \item[(i)] if $x_0\notin\{x_1^*,\ldots,x_k^*\}$,
    \begin{equation}\label{eq:pred_norm_diff}
    \begin{split}
        \Pr(X_{n+1}\in A\mid X^{(n)};\zeta)=&\frac{c}{n}\frac{\int_0^{\infty}u^n\edr^{-c\psi(u)}\tau_1(u)\prod_{m=1}^k\tau_{n_m}(u)\ddr u}{\int_0^{\infty}u^{n-1}\edr^{-c\psi(u)}\prod_{m=1}^{ k}\tau_{n_m}(u)\ddr u}
P_0(A)\\
&+\frac{1}{n}\sum_{l=1}^k \frac{\int_0^{\infty}u^n\edr^{-c\psi(u)}\tau_{n_l+1}(u)\prod_{m\neq l}\tau_{n_m}(u)\ddr u}{ \int_0^{\infty}u^{n-1}\edr^{-c\psi(u)} \prod_{m=1}^k\tau_{n_m}(u)\ddr u}\delta_{x_l^*}(A);
    \end{split}
    \end{equation}
    \item[(ii)] if $x_0=x_j^*$ for some $j=1,\ldots,k$,
\begin{align}\label{eq:pred_norm_ss}
\Pr(X_{n+1}\in A\mid X^{(n)};\zeta)=&\frac{c(1-\zeta)}{n}\frac{\sum_{i=1}^{n_j}c^i\zeta^i\int_0^{\infty}u^n\edr^{-c\psi(u)}\tau_1(u)\xi_{n_j,i}(u)\prod_{m\neq j}\tau_{n_m}(u)\ddr u}{\sum_{i=1}^{n_j}c^i\zeta^i\int_0^{\infty}u^{n-1}\edr^{-c\psi(u)}\xi_{n_j,i}(u)\prod_{m\neq j}\tau_{n_m}(u)\ddr u}
P^*(A)\notag\\
&+\frac{1}{n}\frac{\sum_{i=1}^{n_j+1}c^i\zeta^i \int_0^{\infty}u^n\edr^{-c\psi(u)}\xi_{n_j+1,i}(u)\prod_{m\neq j}\tau_{n_m}(u)\ddr u}{\sum_{i=1}^{n_j}c^i\zeta^i\int_0^{\infty}u^{n-1}\edr^{-c\psi(u)}\xi_{n_j,i}(u)\prod_{m\neq j}\tau_{n_m}(u) \ddr u}\delta_{x_0}(A)\\
&+\frac{1}{n}\sum_{l \neq j}\frac{\sum_{i=1}^{n_j}c^i\zeta^i \int_0^{\infty}u^n\edr^{-c\psi(u)}\tau_{n_l+1}(u)\xi_{n_j,i}(u)\prod_{m\neq l,j}\tau_{n_m}(u)\ddr u}{\sum_{i=1}^{n_j}c^i\zeta^i \int_0^{\infty}u^{n-1}\edr^{-c\psi(u)} \xi_{n_j,i}(u)\prod_{m\neq j}\tau_{n_m}(u)\ddr u}\delta_{x_l^*}(A).\notag
\end{align}
\end{itemize}
\end{theorem}

Note that, if $x_0\not\in\{x_1^*,\ldots,x_k^*\}$, the form of the predictive distribution in \eqref{eq:pred_norm_diff} coincides with that one of the predictive distribution of a hNRMI with diffuse base measure (see \citealp[Corollary~1]{Jam06}).

We complete this section by studying the distribution of $N_0^{(n)}$, the number of elements, in an  exchangeable sample $X^{(n)}$ such that $X_i\mid \tilde P\simiid \tilde P$, that coincide with $x_0$. 
\begin{theorem}\label{thm:N0_norm_ss}
Let $X_i\mid \tilde P\simiid \tilde P$, for $i=1,\ldots,n$ and for any $n\ge 1$. If $\tilde P\sim\mbox{\rm hNRMI}(\rho;c,P_0)$, with $P_0$ as in \eqref{eq:base_inner}, then the distribution of $N_0^{(n)}$, parametrized by $\zeta\in(0,1)$, is given, for any $j\in\{0,1,\ldots,n\}$, by
\begin{equation}
\begin{split}
    \Pr&(N_0^{(n)}=j;\zeta)=\frac{1}{\Gamma(n)}\sum_{k=1}^{n-j+1-\delta_{0j}} k^{1-\delta_{0j}} c^{k+\delta_{0j}-1} (1-\zeta)^{k+\delta_{0j}-1}\\
    &\times \sum_{i=0}^j c^i\zeta^i \int_0^\infty u^{n-1}\edr^{-c\psi(u)}\xi_{j,i}(u)\frac{1}{(k+\delta_{0j}-1)!}\sum_{\bm{n}\in \mathcal{N}_{k+\delta_{0j}-1}^{(n-j)}}\binom{n}{n_1\cdots n_{k+\delta_{0j}-1}}\prod_{m=1}^{k+\delta_{0j}-1}\tau_{n_m}(u)\ddr u,
\end{split}\label{eq:N0_norm_ss}
\end{equation}
where $\mathcal{N}_k^{(n)}=\{\bm{n}=(n_1,\ldots,n_k)\,:\,n_m \in \mathds{N} \text{ and } \sum_{m=1}^k n_m=n\}$, and  $\delta_{0j}=1$ if $j=0$ and $\delta_{0j}=0$ otherwise.
\end{theorem}

The expressions obtained for the case of general hNRMIs with spike and slab base measure showcase how the techniques developed in the paper work in great generality for a large class of inner spike and slab models. In order to make such expressions amenable of direct application, special cases of the class of hNRMIs must be considered, as illustrated in the next section.

\section{Notable special cases}\label{sec:stable}
Here we specialize the results of Section \ref{sec:main} to two popular special cases, namely the $\sigma$-stable hNRMI and the N-IG process. This allows us to display that the expressions we obtained for the EPPF, the predictive distributions, and the distribution of $N_0^{(n)}$, for the general case of hNRMIs, reduce to tractable expressions when special cases within the same family are considered. Moreover, the $\sigma$-stable hNRMI is the only random probability measure, along with the DP, which is at the same time a special case of the Pitman-Yor process and an element of the class of NRMIs. The last observation allows us to link the results of this paper with those in \citet{Can17}. 

\subsection{$\sigma$-stable hNRMI}
Assume that $\tilde P$ is a $\sigma$-stable hNRMI, with $\sigma\in(0,1)$, with spike and slab base measure. In other terms, $\tilde P\sim \text{hNRMI}(\rho;1,P_0)$ with $\rho(s)=\sigma s^{-1-\sigma} /\Gamma(1-\sigma)$, $c$ set equal to 1, and $P_0$ defined as in \eqref{eq:base_inner}. Such choice implies that $\psi(u)=u^\sigma$, $\tau_q(u)=u^{-q+\sigma} \sigma (1-\sigma)_{q-1}$, and $\xi_{n,i}(u)=u^{\sigma i-n}\mathscr{C}(n,i;\sigma)$, where $\mathscr{C}(n,i;\sigma)=\frac{1}{i!}\sum_{r=0}^i(-1)^r \binom{i}{r}(-r\sigma)_{n}$ is the generalized factorial coefficient \citep{Cha05}, and $(a)_n=\Gamma(a+n)/\Gamma(a)$. These expressions for $\psi(u)$, $\tau_n(u)$ and $\xi_{n,i}(u)$, when plugged into \eqref{eq:eppf:nrmi}, 
\eqref{eq:pred_norm_diff}, \eqref{eq:pred_norm_ss} and \eqref{eq:N0_norm_ss}, provide the EPPF, the prediction rule, and the distribution of $N_0^{(n)}$, implied by the $\sigma$-stable hNRMI $\tilde P$. For the sake of compactness we introduce the numbers $\varphi_{m,q}$, defined, for any 
positive integers $m$ and $q$, as
$$\varphi_{m,q}(\zeta)=\sum_{i=1}^m \zeta^i \Gamma(q+i-1)\mathscr{C}(m,i;\sigma).$$ 
In the following, for the sake of simplicity, we will omit the dependence on $\zeta$ and write $\varphi_{m,q}$ instead of $\varphi_{m,q}(\zeta)$.
The EPPF can then be written as
	\begin{equation*}\label{eq:eppf:stableNRMI}
	\Pi_k^{(n)}(n_1,\ldots,n_k;\zeta)=\frac{\sigma^{k-2}(1-\zeta)^{k-1}}{\Gamma(n)}\prod_{m=1}^k (1-\sigma)_{n_m-1}\left(\sigma(1-\zeta) \Gamma(k)+\sum_{l=1}^{k}\frac{\varphi_{n_l,k}}{(1-\sigma)_{n_l-1}}\right).
	\end{equation*}
	The predictive distributions of Theorem \ref{thm:predictive} reduce to:
	\begin{itemize}
    \item[(i)] if $x_0\notin\{x_1^*,\ldots,x_k^*\}$,
    \begin{equation}\label{eq:pred_stable_diff}
        \Pr(X_{n+1}\in A\mid X^{(n)};\zeta)=\frac{k\sigma}{n}P_0(A)+\frac{1}{n}\sum_{l=1}^k (n_{l}-\sigma)\delta_{x_l^*}(A);
    \end{equation}
    \item[(ii)] if $x_0=x_j^*$ for some $j=1,\ldots,k$,
\begin{align}\label{eq:pred_stable_ss}
\Pr(X_{n+1}\in A\mid X^{(n)};\zeta)=&\frac{(1-\zeta)\sigma}{n}\frac{\varphi_{n_j,k+1}}{\varphi_{n_j,k}}
P^*(A)+\frac{1}{n}\frac{\varphi_{n_j+1,k}}{\varphi_{n_j,k}}\delta_{x_0}(A)+\frac{1}{n}\sum_{l \neq j}(n_l-\sigma)\delta_{x_l^*}(A).
\end{align}
\end{itemize}
We observe that \eqref{eq:pred_stable_diff} has the exact same structure of the predictive distribution of a $\sigma$-stable hNRMI with diffuse base measure. The same does not hold for the predictive distribution in \eqref{eq:pred_stable_ss}, which accounts for the fact that $X_{n+1}$ can coincide with $x_0$ either when it is drawn from the spike component of the base measure $P_0$ or when it is tied to $x_j^*$. It is interesting to observe that, by setting $A=\X$ in \eqref{eq:pred_stable_ss}, the next triangular identity is obtained
\begin{equation*}
    \varphi_{m,q}(m+(q-1)\sigma)=\varphi_{m+1,q}+(1-\zeta)\sigma \varphi_{m,q+1},
\end{equation*}
for any positive integers $m$ and $q$. This, combined with the fact that  $\varphi_{1,q}=\zeta \sigma \Gamma(q)$, for any positive integer $q$, suggests a recursive strategy for efficiently evaluate the numbers $\varphi_{m,q}$, for the desired values of $m$ and $q$.\\
Finally, the distribution of $N_0^{(n)}$ provided in Theorem \ref{thm:N0_norm_ss} becomes
\begin{equation}\label{eq:N0_sigma_ss}
    \Pr(N_0^{(n)}=j;\zeta)=\frac{1}{\sigma \Gamma(n)}\sum_{k=1}^{n-j+1-\delta_{0j}} k^{1-\delta_{0j}} (1-\zeta)^{k+\delta_{0j}-1} \mathscr{C}(n-j,k+\delta_{0j}-1;\sigma) \sum_{i=0}^j \zeta^i \Gamma(k+\delta_{0j}+i-1)\mathscr{C}(j,i;\sigma),
\end{equation}
for any $j\in\{0,1,\ldots,n\}$.

\subsection{Normalized inverse Gaussian process}
Assume that $\tilde P$ is a N-IG process, with $\tau>0$ and spike and slab base measure. In other terms, $\tilde P\sim \text{hNRMI}(\rho;c,P_0)$ with $\rho(s)=s^{-3/2}\edr^{-\tau s}/(2\sqrt{\pi})$, $c>0$ and $P_0$ defined as in \eqref{eq:base_inner}. 
Such choice implies that $\psi(u)=\sqrt{\tau+u}-\sqrt{\tau}$, $\tau_q(u)=(\tau+u)^{1/2-q}\Gamma(q-1/2)/(2\sqrt{\pi})$, and $\xi_{n,i}(u)=(\tau+u)^{i/2-n}\mathscr{C}(n,i;1/2)$. These expressions for $\psi(u)$, $\tau_n(u)$ and $\xi_{n,i}(u)$, when plugged into \eqref{eq:eppf:nrmi}, 
\eqref{eq:pred_norm_diff}, \eqref{eq:pred_norm_ss}, and \eqref{eq:N0_norm_ss}, provide the EPPF, the prediction rule, and the distribution of $N_0^{(n)}$, implied by a normalized inverse Gaussian process. For the sake of compactness, we introduce the notation $\beta=c\sqrt{\tau}$ and define the numbers $\varrho_{m,q}^{(n)}$ as follows. For any $n\geq q$ positive integers,
      $$\varrho_{0,q}^{(n)}=\sum_{r=0}^{n-1}\binom{n-1}{r}(-1)^r\beta^{2r}\Gamma(q-2r,\beta),$$
      where $\Gamma(n,z)=\int_z^{\infty}x^{n-1}\edr^x\ddr x$ denotes the upper incomplete gamma function; moreover, for any $m$, $q$ and $n$ positive integers such that $n\geq q$ and $m\leq n-(q-1)$,
      $$\varrho_{m,q}^{(n)}(\zeta)=\sum_{i=1}^{m}(2\zeta)^i\frac{\Gamma(2m-i)}{\Gamma(m+1-i)\Gamma(i)}\varrho_{0,q+i-1}^{(n)}.$$
With the aim of simplifying the notation, in the following we will omit the dependence on $\zeta$ and write $\varrho_{m,q}^{(n)}$ instead of $\varrho_{m,q}^{(n)}(\zeta)$.
The EPPF can then be written as	\begin{equation*}\label{eq:eppf:N-IG}
	\Pi_k^{(n)}(n_1,n_2,\ldots,n_k;\zeta)=\frac{\edr^{\beta}(1-\zeta)^{k-1}}{2^{k-1}\pi^{k/2} \Gamma(n)}\prod_{m=1}^k\Gamma(n_m-1/2)\left( (1-\zeta)\varrho_{0,k}^{(n)}+\sqrt{\pi}\sum_{l=1}^k \frac{2^{1-2n_l}}{\Gamma(n_l-1/2)}\varrho_{n_l,k}^{(n)}\right).
	\end{equation*}
	The predictive distribution of Theorem \ref{thm:predictive} simplifies to:
	\begin{itemize}
    \item[(i)] if $x_0\notin\{x_1^*,\ldots,x_k^*\}$,
    \begin{equation}\label{eq:pred_NIG_diff}
        \Pr(X_{n+1}\in A\mid X^{(n)};\zeta)=\frac{1}{2n} \frac{\varrho_{0,k+1}^{(n+1)}}{\varrho_{0,k}^{(n)}}P_0(A)+\frac{1}{n}\frac{\varrho_{0,k}^{(n+1)}}{\varrho_{0,k}^{(n)}}\sum_{l=1}^k (n_{l}-1/2)\delta_{x_l^*}(A);
    \end{equation}
    \item[(ii)] if $x_0=x_j^*$ for some $j=1,\ldots,k$,
\begin{equation}\label{eq:pred_NIG_ss}
\Pr(X_{n+1}\in A\mid X^{(n)};\zeta)=\frac{(1-\zeta)}{2n}\frac{\varrho_{n_j,k+1}^{(n+1)}}{\varrho_{n_j,k}^{(n)}}
P^*(A)+\frac{1}{4n}\frac{\varrho_{n_j+1,k}^{(n+1)}}{\varrho_{n_j,k}^{(n)}}\delta_{x_0}(A)+\frac{\varrho_{n_j,k}^{(n+1)}}{\varrho_{n_j,k}^{(n)}}\sum_{l \neq j}\frac{(n_l-1/2)}{n}\delta_{x_l^*}(A).
\end{equation}
\end{itemize}
We observe that \eqref{eq:pred_NIG_diff} displays the standard structure of the predictive distribution of a N-IG process with diffuse base measure \citep[see][]{Lij05}. The same does not hold for the predictive distribution in \eqref{eq:pred_NIG_ss}, which accounts for the fact that $X_{n+1}$ can coincide with $x_0$ either when it is drawn from the spike component of the base measure $P_0$ or when it is tied to $x_j^*$.\\ 
It is interesting to observe that, by setting $A=\X$ in \eqref{eq:pred_NIG_diff}, the next triangular identity is obtained
\begin{equation*}
2n\varrho_{0,q}^{(n)}=\varrho_{0,q+1}^{(n+1)}+(2n-q)\varrho_{0,q}^{(n+1)},
\end{equation*}
for any positive integers $n\geq q$. In a similar fashion, by setting $A=\X$ in \eqref{eq:pred_NIG_ss}, one can show that the identity
\begin{equation*}
    4n \varrho_{m,q}^{(n)}=2(1-\zeta)\varrho_{m,q+1}^{(n+1)}+\varrho_{m+1,q}^{(n+1)}+2(2(n-m)-q+1)\varrho_{m,q}^{(n+1)}
\end{equation*}
holds true for any positive integers $n,m,q$ such that $q\leq n$ and $1\leq m\leq n-q+1$.

Finally, the distribution of $N_0^{(n)}$ provided in Theorem \ref{thm:N0_norm_ss} becomes
\begin{equation}
    \Pr(N_0^{(n)}=j;\zeta)=\frac{2\edr^\beta}{\Gamma(n)}\sum_{k=1}^{n-j+1-\delta_{0j}} k^{1-\delta_{0j}} (1-\zeta)^{k+\delta_{0j}-1} \mathscr{C}(n-j,k+\delta_{0j}-1;1/2) \sum_{i=0}^j \zeta^i \mathscr{C}(j,i;1/2)\varrho_{0,k+i-1}^{(n)},
\end{equation}
for any $j\in\{0,1,\ldots,n\}$.

\section{Generalized P\'olya urn scheme and numerical study}\label{sec:algo}
In this section we describe the steps of a generalized P\'olya urn scheme which, under the assumption that $X_i\mid \tilde P \simiid \tilde P$, with $\tilde P$ distributed as a hNRMI with spike and slab base measure \eqref{eq:base_inner}, allows us to generate a sample $X_{n+1},\ldots,X_{n+m}$, conditionally on the observation of $(x_1,\ldots,x_n)$. The steps are summarized 
in Algorithm \ref{alg:PUS}.

\begin{algorithm}[h!]
\SetAlgoLined
Set $r=1$\;
\While{$r\leq m$}{
let $k$ be the number of distinct values $(x_1^*,\ldots,x_k^*)$ in $(x_1,\ldots,x_n,x_{n+1},\ldots,x_{n+r-1})$, and let  $(n_1,\ldots,n_k)$ be the corresponding frequencies\;
  \eIf{$x_0\neq x_j^*$ for every $j\in\{1,\ldots,k\}$}{
  generate a realization $x_{n+r}$ of $X_{n+r}$, conditionally on the observation of $X^{(n+r-1)}=(x_{1},x_{2},\ldots,x_{n+r-1})$, from the predictive distribution \eqref{eq:pred_norm_diff}\;
   }{
   generate a realization $x_{n+r}$ of $X_{n+r}$, conditionally on the observation of $X^{(n+r-1)}=(x_{1},x_{2},\ldots,x_{n+r-1})$, from the predictive distribution \eqref{eq:pred_norm_ss}\;
 }
 increase $r$ by 1\;
 }
 \KwRet{$(x_{n+1},\ldots,x_{n+m})$ 
    }
 \caption{Generalised P\'olya urn scheme for hNRMI inner spike and slab model.}\label{alg:PUS}
\end{algorithm}

While such a scheme can be adopted, in line of principle, for any specification of the L\'evy intensity $\rho$, the results displayed in Section \ref{sec:stable} allow for a direct implementation of this scheme for $\tilde P$ distributed as a $\sigma$-stable hNRMI or a N-IG process.

Next we implement the scheme we devised, in order to complement with a simulation study the findings of Proposition \ref{prop:var} on the different features of inner and outer spike and slab models. For the purpose of illustration, we focus on the case of $\sigma$-stable hNRMIs. In this study we investigate the variability of both the prior and the posterior distributions obtained by considering inner and outer spike and slab models, $\tilde P$ and $\tilde Q$ respectively, for different specifications of the parameters $\zeta$ and $\sigma$.

\subsection{Inner vs outer models a priori}
\label{sec:study1} 
\label{sec:priorstudy}

When analysing data, prior predictive checks represent a very effective, yet easy to implement, strategy to get a handle on the implications of a specific prior choice on functionals of interest. With this in mind, we study and compare the prior behaviour of inner and outer spike and slab models by generating, under a range of settings, realizations $(x_{1},\ldots,x_{m})$ of the random vector $(X_1,\ldots,X_m)$, modelled with either $\tilde P$ or $\tilde Q$. We then consider the corresponding empirical cumulative distribution function $F_m(x)=m^{-1}\sum_{i=1}^m \mathds{1}_{(-\infty,x_i]}(x)$, realization of the random cumulative distribution function $\tilde{F}_m(x)=m^{-1}\sum_{i=1}^m \mathds{1}_{(-\infty,X_i]}(x)$. We replicate the same experiment 100\,000 times, which allows us to quantify, for any $x' \in\R$, the uncertainty associated with the random variable $\tilde{F}_m(x')$, by means of a Monte Carlo estimate of a quantile-based $95\%$ credible interval. While both inner and outer spike and slab model are such that $\E[\tilde F_m(x)]=F_0(x)$, where $F_0$ denotes the cumulative distribution function corresponding to the spike and slab base measure \eqref{eq:base_inner}, our simulation compares the associated uncertainties and investigates the role played by the parameters $\zeta$ and $\sigma$.

\begin{figure}[h!]
\centering
\includegraphics[width=0.25\textwidth,bb=260 10 990 700,clip=true]{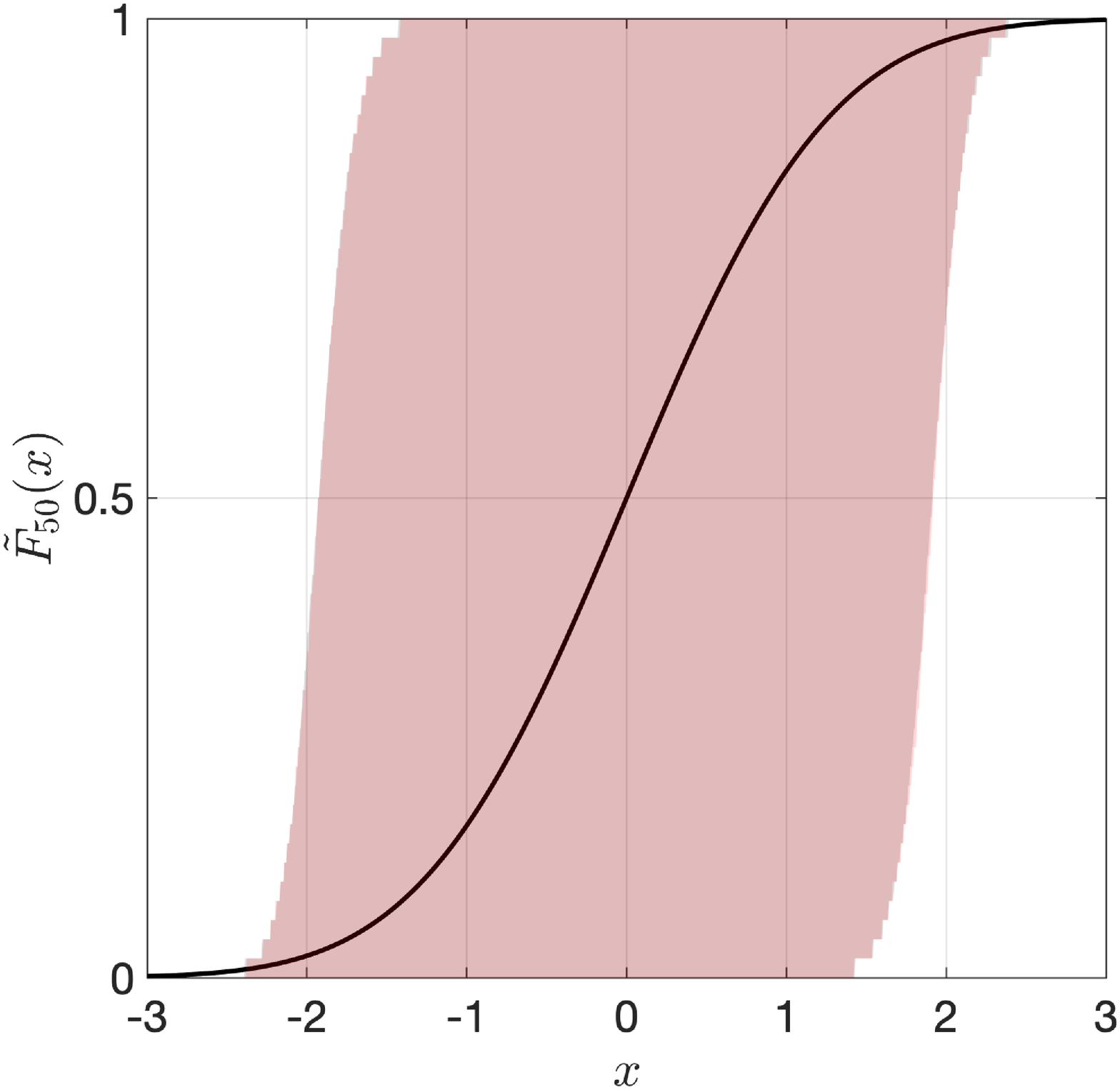}
\hspace{-0.1cm}%
\includegraphics[width=0.25\textwidth,bb=260 10 990 700,clip=true]{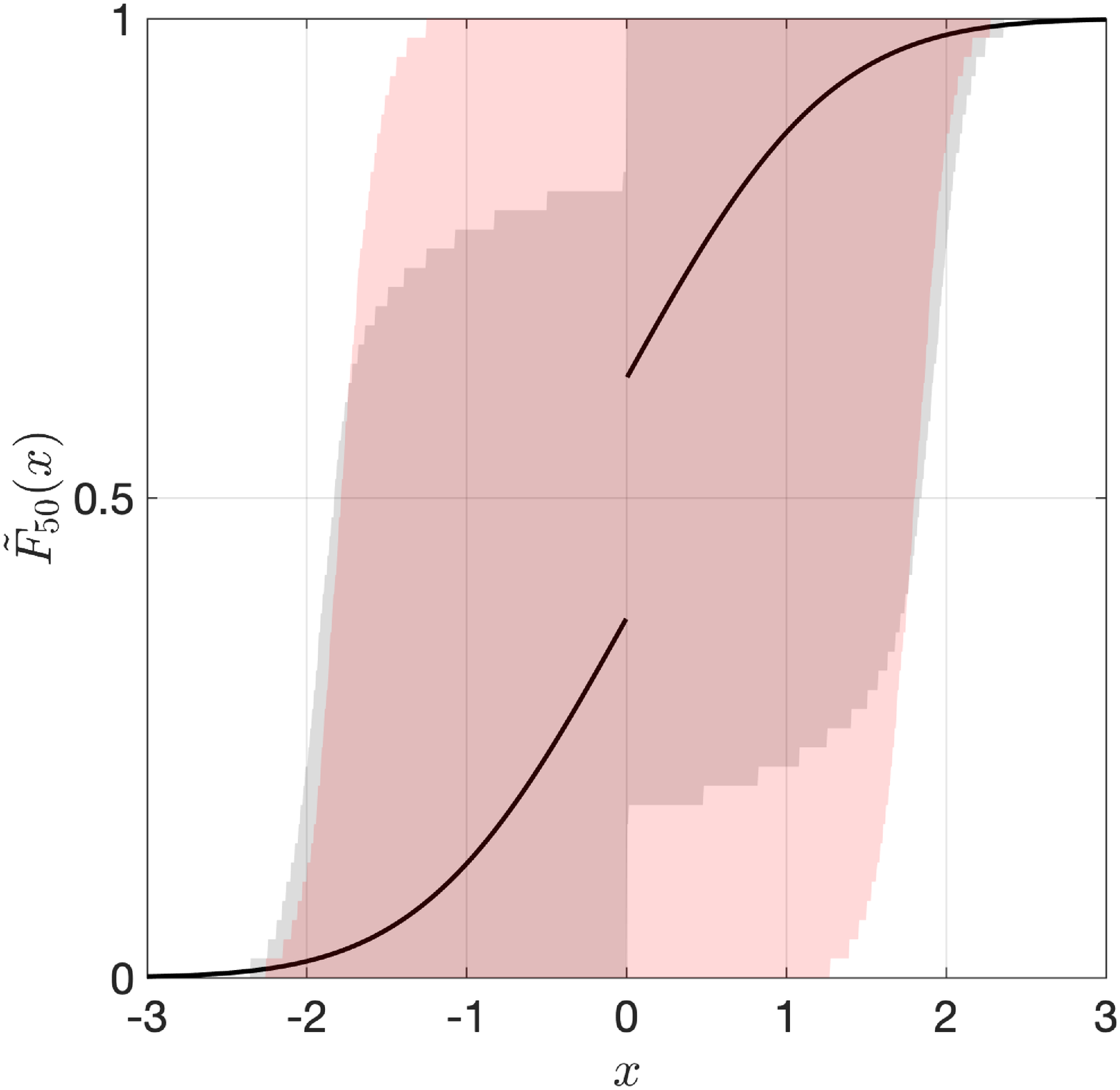}
\hspace{-0.1cm}%
\includegraphics[width=0.25\textwidth,bb=260 10 990 700,clip=true]{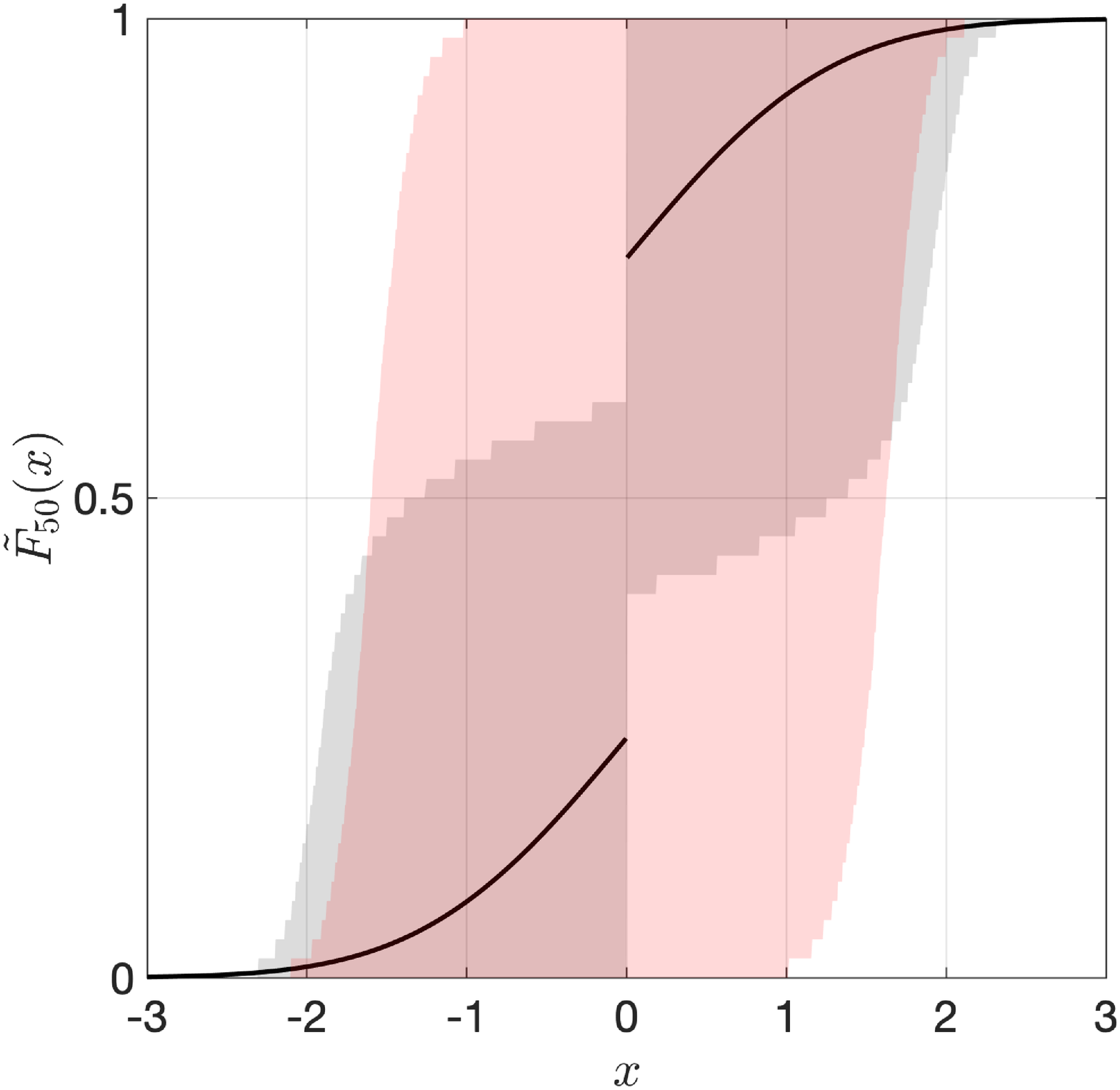}
\hspace{-0.1cm}%
\includegraphics[width=0.25\textwidth,bb=260 10 990 700,clip=true]{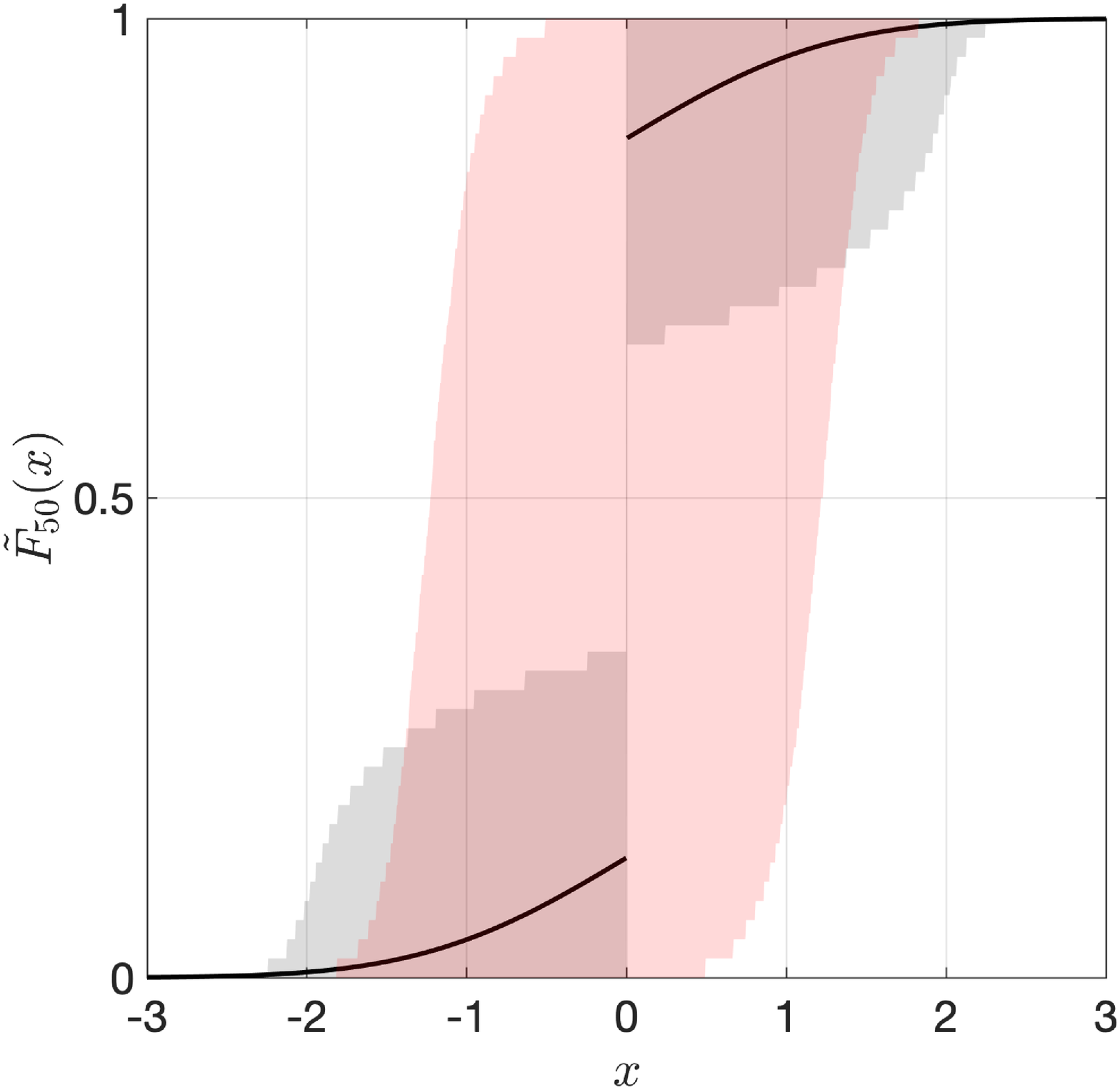}\\[10pt] \includegraphics[width=0.25\textwidth,bb=260 10 990 700,clip=true]{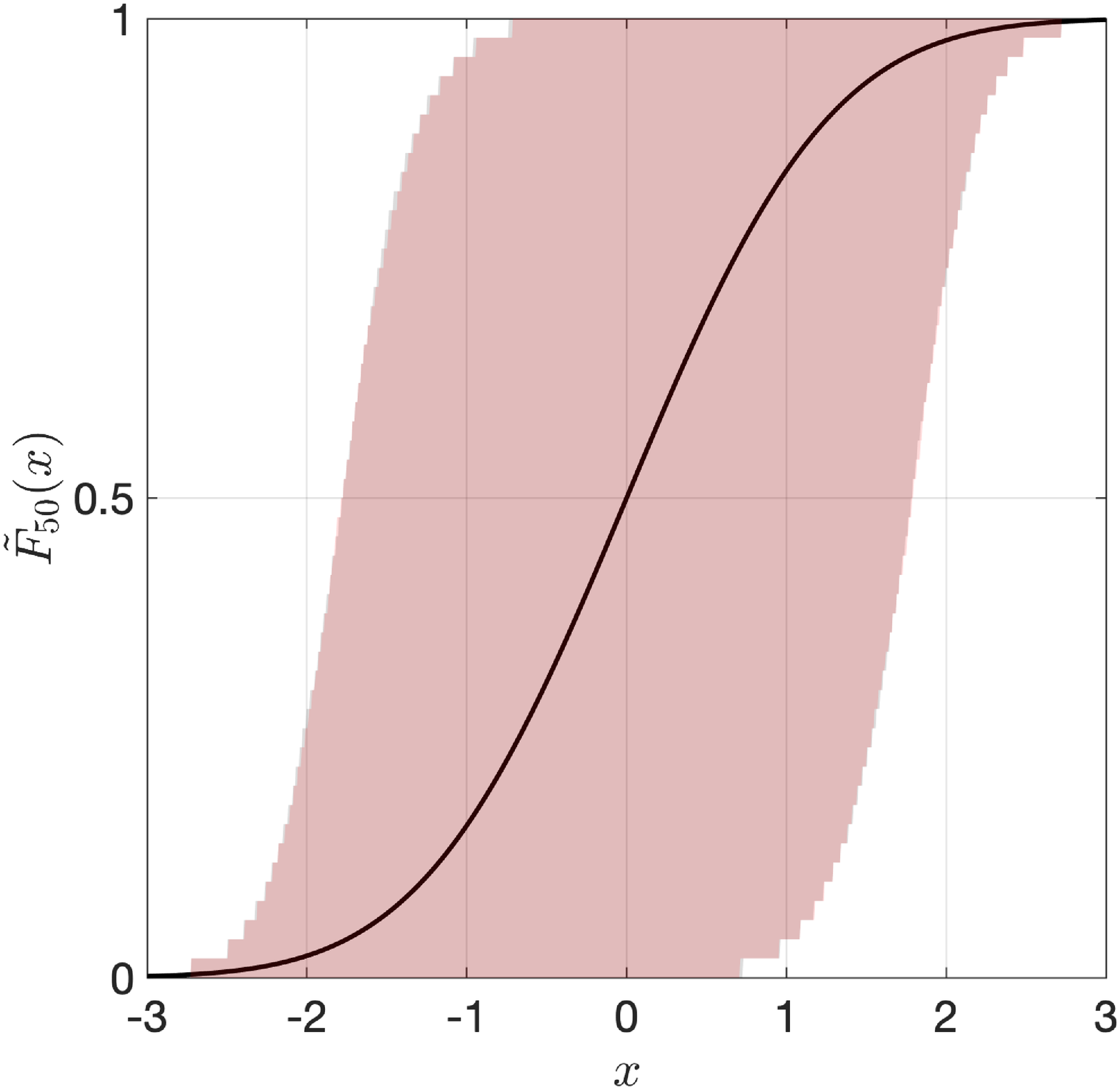}
\hspace{-0.1cm}%
\includegraphics[width=0.25\textwidth,bb=260 10 990 700,clip=true]{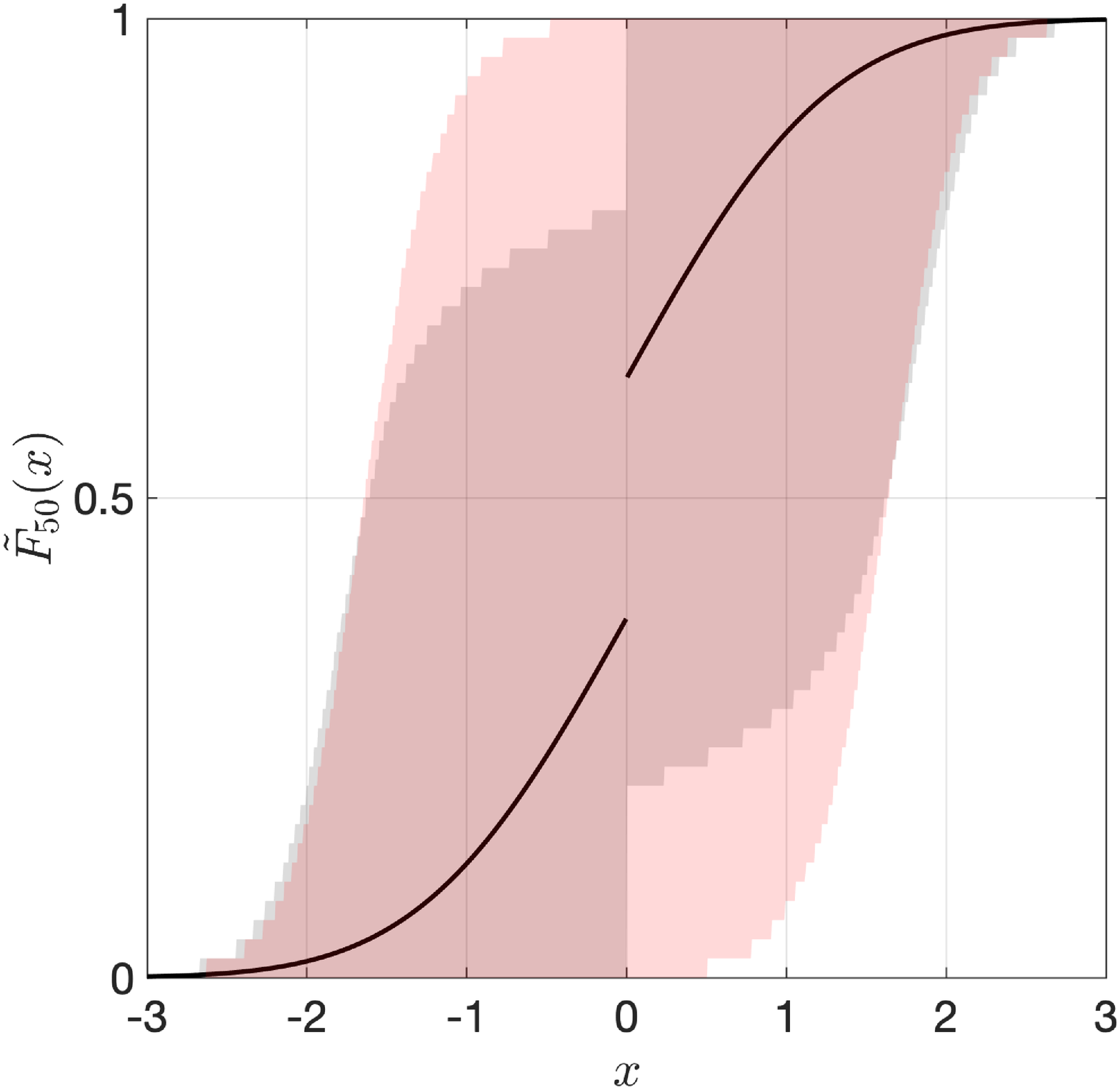}
\hspace{-0.1cm}%
\includegraphics[width=0.25\textwidth,bb=260 10 990 700,clip=true]{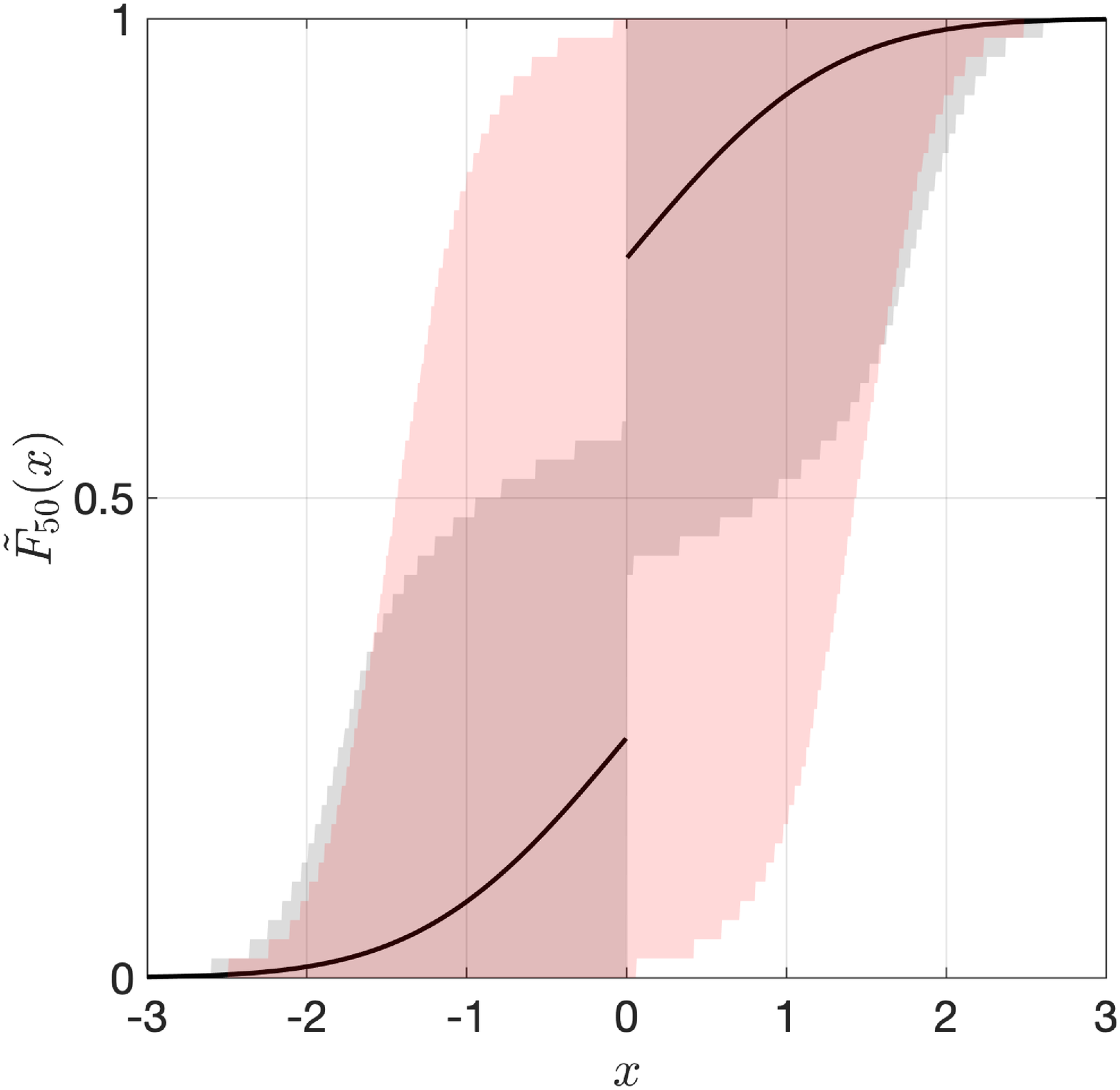}
\hspace{-0.1cm}%
\includegraphics[width=0.25\textwidth,bb=260 10 990 700,clip=true]{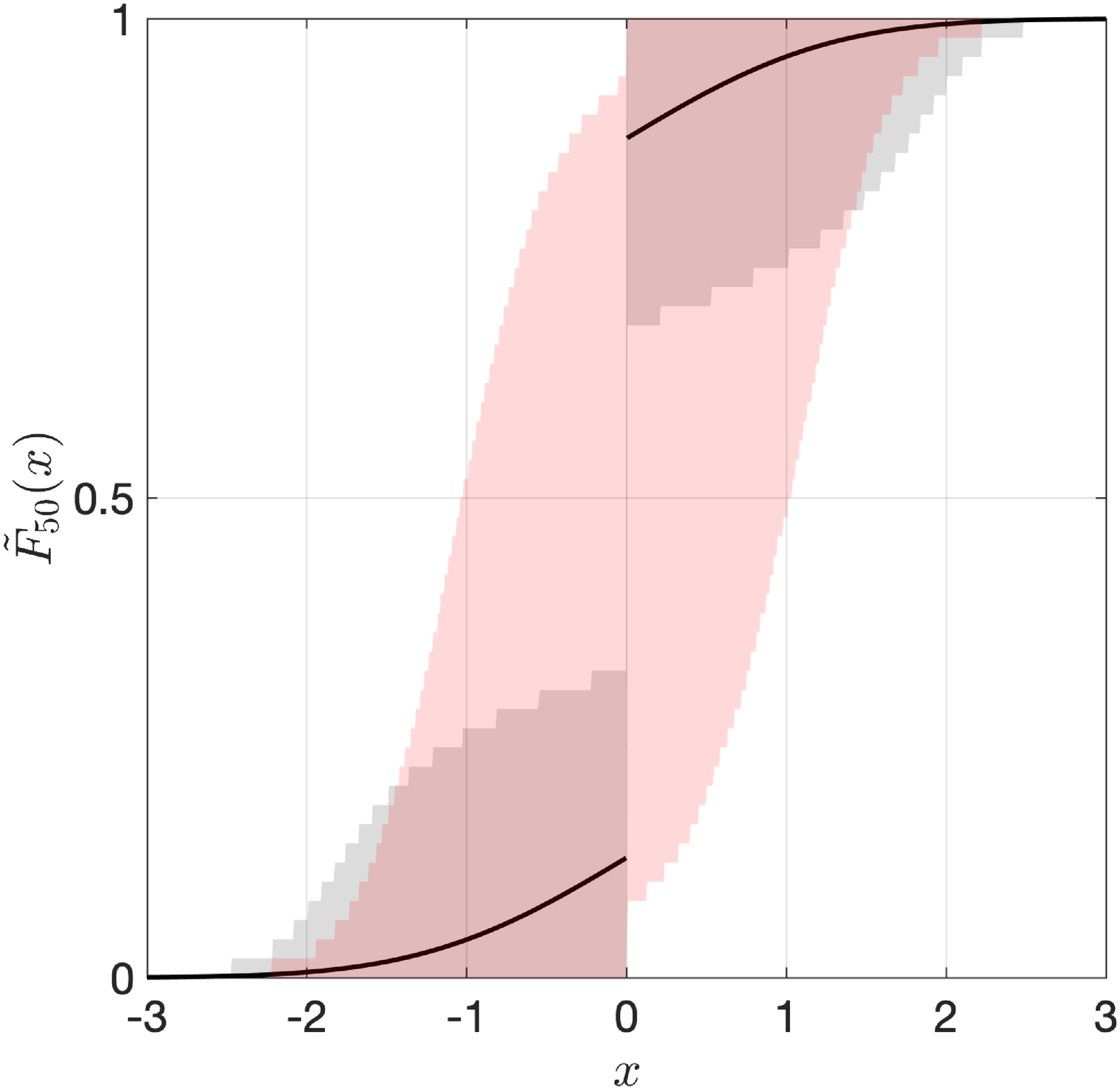}\\[10pt]
\includegraphics[width=0.25\textwidth,bb=260 10 990 700,clip=true]{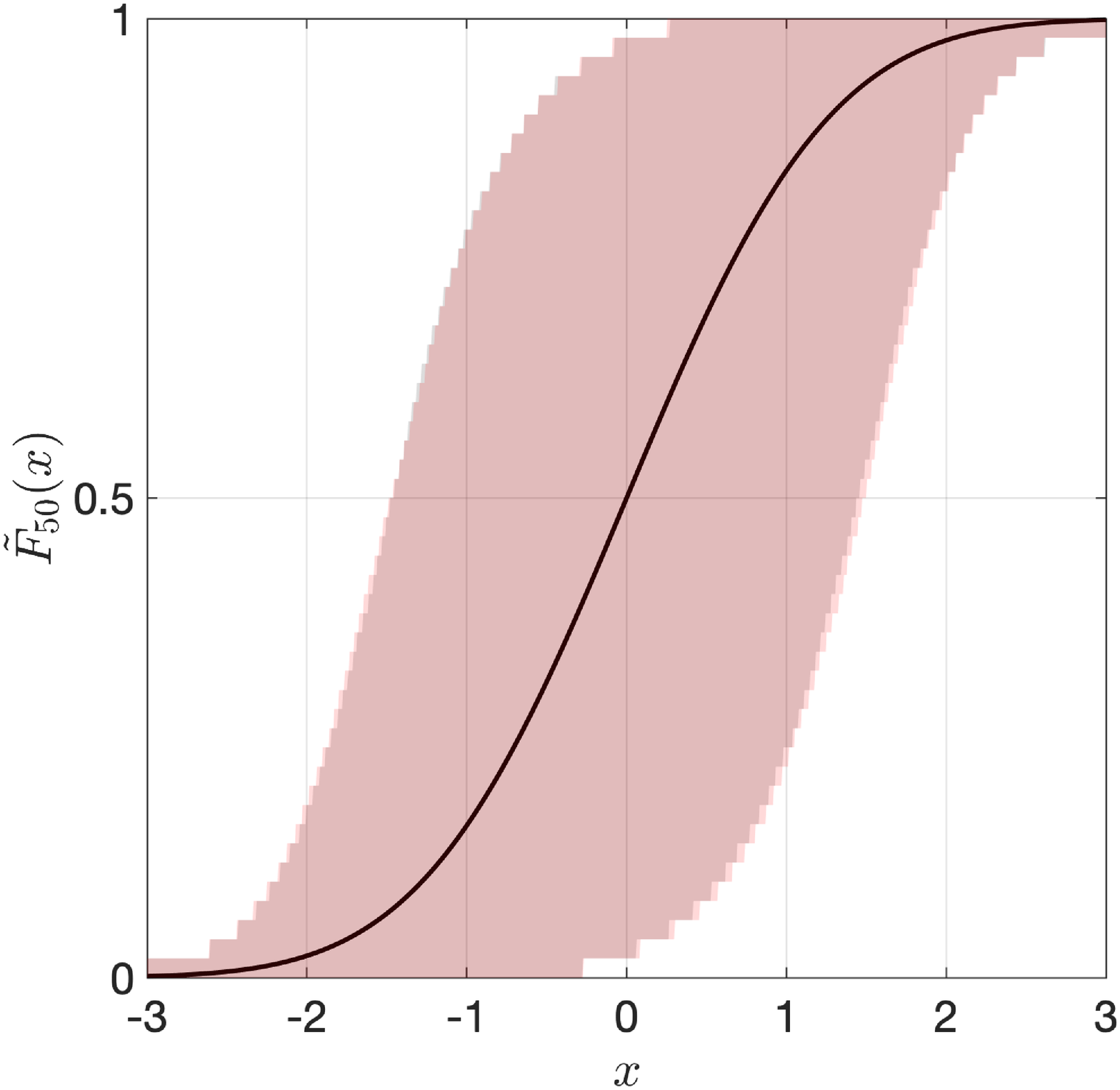}
\hspace{-0.1cm}%
\includegraphics[width=0.25\textwidth,bb=260 10 990 700,clip=true]{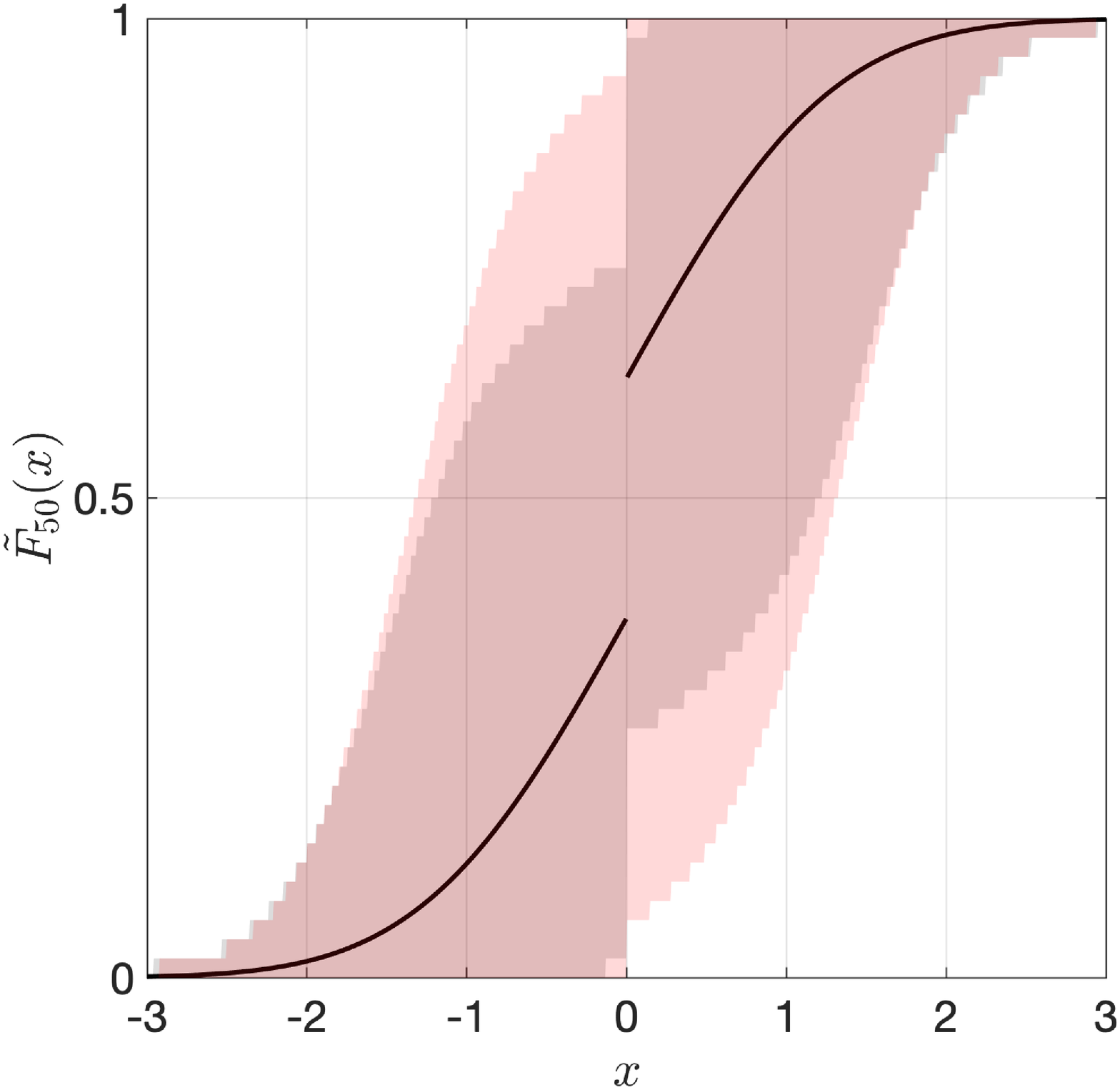}
\hspace{-0.1cm}%
\includegraphics[width=0.25\textwidth,bb=260 10 990 700,clip=true]{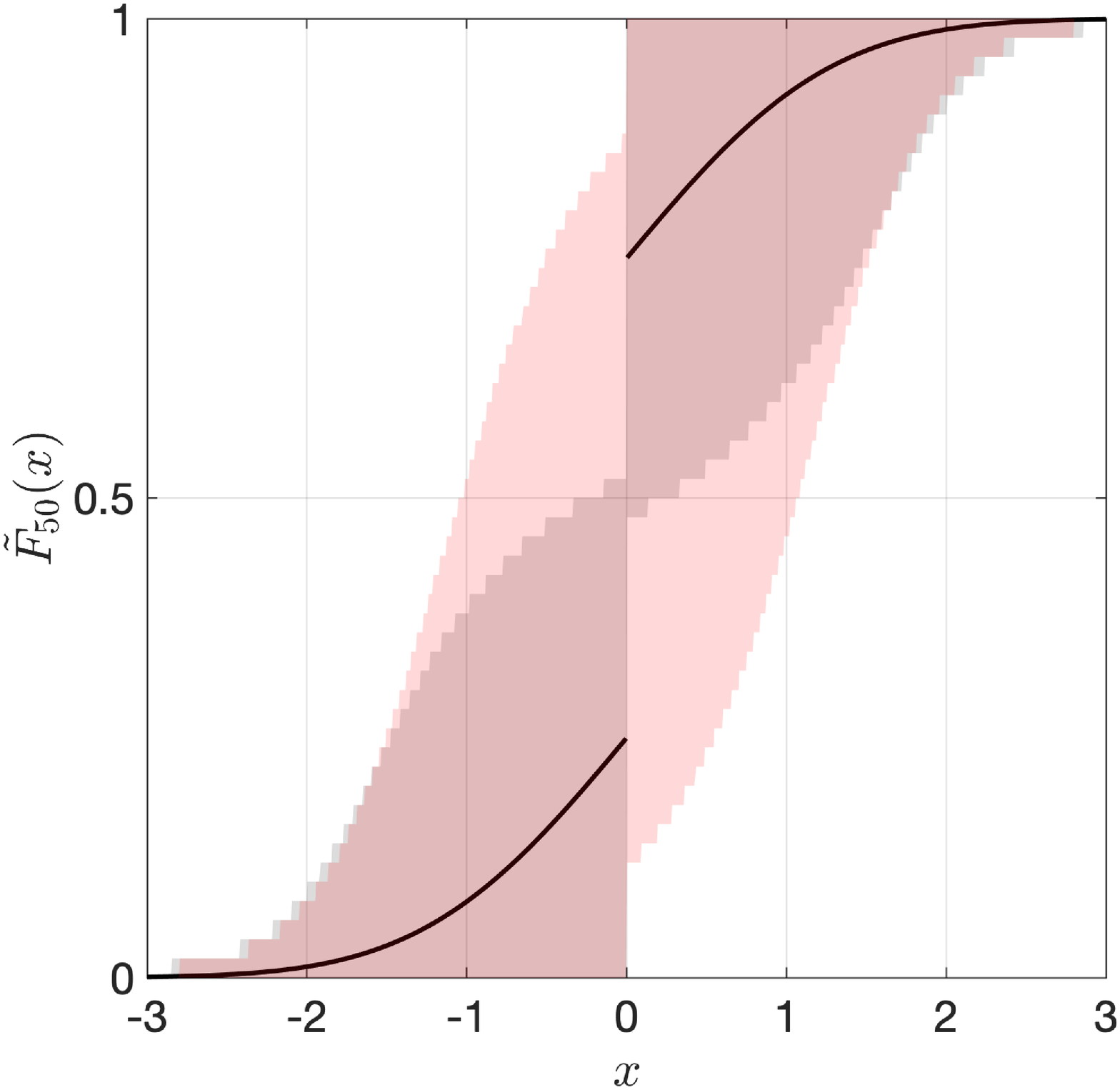}
\hspace{-0.1cm}%
\includegraphics[width=0.25\textwidth,bb=260 10 990 700,clip=true]{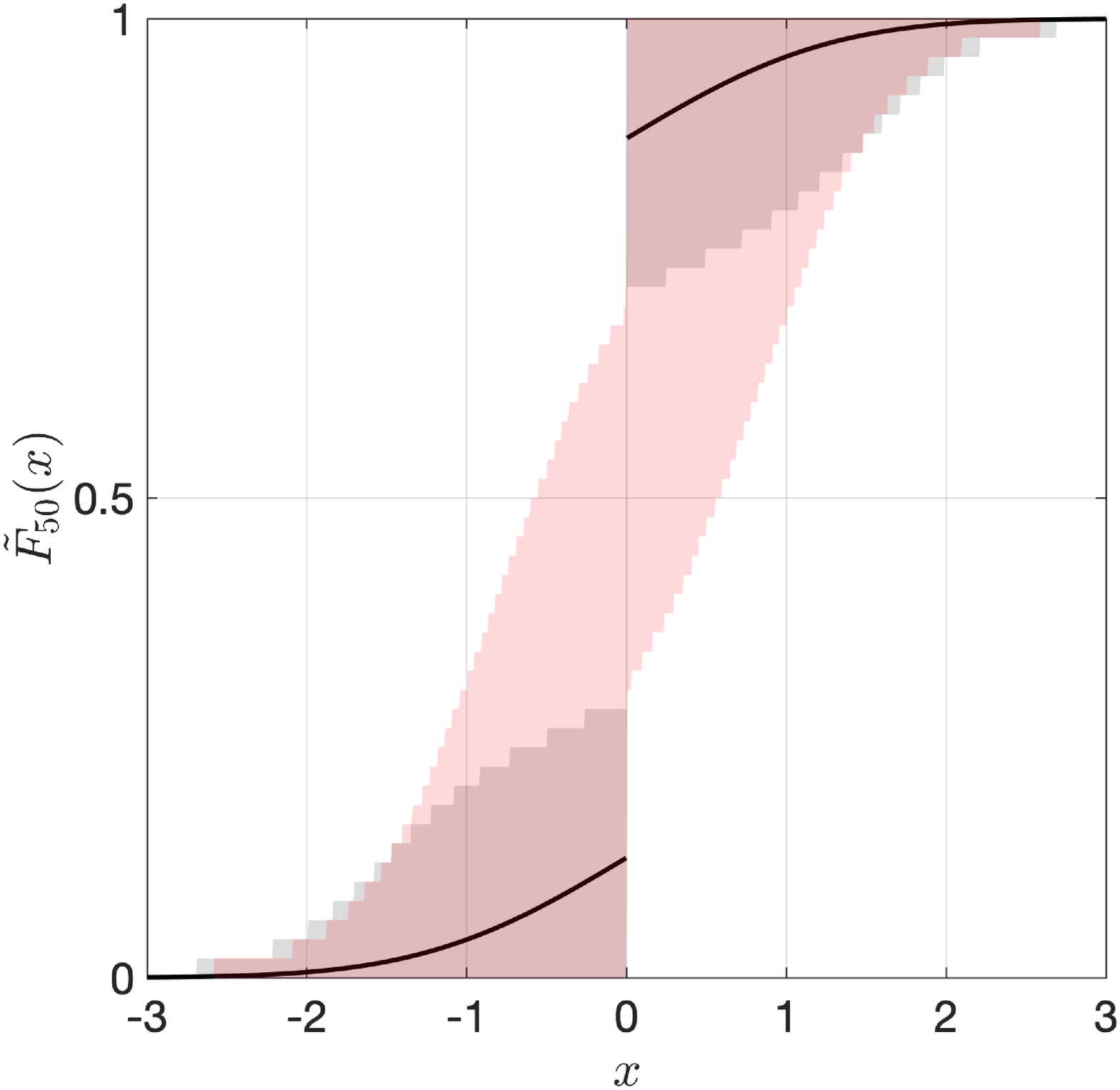}
\caption{Comparison of the estimated quantile-based $95\%$ credible bands for the prior distributions of $\tilde F_{50}$ induced by inner (red shades) and outer (grey shades) spike and slab models for $\sigma$-stable hNRMIs, along with its expected value (black curves) $\mathds{E}[\tilde F_{50}]=F_0$. Different rows refer to different values of $\sigma\in\{0.25,0.5,0.75\}$ (from top to bottom), different columns refer to different values of $\zeta\in\{0,0.25,0.5,0.75\}$ (from left to right).}\label{fig:only}
\end{figure}

Let $P^*$ be a standard normal distribution, $x_0=0$, and set $m=50$. We consider $\zeta\in\{0,0.25,0.5,0.75\}$ and $\sigma\in\{0.25,0.5,0.75\}$. For each of the resulting 12 combinations of parameter values, Figure \ref{fig:only} displays the expected value $F_0$ of $\tilde F_{50}$, as well as the estimated quantile-based $95\%$ credible bands for the inner and the outer spike and slab models. First we observe that, when $\zeta=0$, the two models coincide and boil down to a $\sigma$-stable hNRMI with diffuse base measure $P^*$. On the other hand, it can be appreciated that, when $\zeta>0$, the variability of $\tilde F_{50}$ appears larger for the inner specification of the spike and slab model, such difference being more evident for larger values of $\zeta$. The same behaviour can be appreciated across different values of the parameter $\sigma$, with larger values of $\sigma$ leading to an overall smaller prior variability for both models.
The structural difference between $\tilde P$ and $\tilde Q$ is further highlighted by looking at the distribution of two functionals of $\tilde F_{50}$, namely the mean $M_{50}=\sum_{i=1}^{50}X_i/50$ and the median $\text{med}_{50}=\text{median}(X_1,\ldots,X_{50})$. Table \ref{tab:only} reports the length of the estimated $95\%$ quantile-based credible intervals for $M_{50}$ and for $\text{med}_{50}$, for inner and outer spike and slab models for the same combination of values for $\zeta$ and $\sigma$ considered above. While all the estimated intervals are approximately centered at $x_0=0$, their lengths vary considerably. For both $M_{50}$ and $\text{med}_{50}$, it can be appreciated that, for a fixed $\zeta$, larger values of $\sigma$ correspond to smaller intervals; similarly, for a fixed $\sigma$, larger values of $\zeta$ display smaller intervals. Intervals corresponding to inner and outer spike and slab models coincide when $\zeta=0$, while the latter tend to be smaller for larger values of $\zeta$. Such difference is particularly sizeable when $\zeta=0.75$, even more so if the median is considered: while the length of the estimated interval for the distribution of $\text{med}_{50}$ induced by $\tilde P$ is equal to 2.44, 2.04 and 1.14 when $\sigma$ takes value 0.25, 0.5 and 0.75, respectively, the corresponding intervals for the distribution of $\text{med}_{50}$ induced by $\tilde Q$ appear degenerate at $x_0=0$, which implies that $\Pr(\text{median}(X_1,\ldots,X_{50})=0)>0.95$.

\begin{table}[h!]
\makegapedcells
\centering
\begin{tabular}{ c c | c c c c | c c c c | c c c c }
  \multicolumn{2}{c|}{$\sigma$} & \multicolumn{4}{c|}{$0.25$} & \multicolumn{4}{c|}{$0.5$} & \multicolumn{4}{c}{$0.75$} \\
  \hline
  \multicolumn{2}{c|}{$\zeta$} & $0$ &  $0.25$ &  $0.5$ &  $0.75$ & $0$ &  $0.25$ &  $0.5$ &  $0.75$ & $0$ &  $0.25$ &  $0.5$ &  $0.75$\\
  \hline
  \multirow{2}{*}[-0.3em]{\rotatebox[origin=c]{90}{$M_{50}$}} & inner & 3.48 &    3.22 &    2.86 &    2.14 &    2.93 &    2.67 &    2.31 &    1.65 &    2.20 &    1.95 &    1.59 &    1.10  \\
  & outer & 3.48 &    2.60 &    1.77 &    0.92 &    2.94 &    2.21 &    1.49 &    0.78 &    2.17 &   1.64 &    1.12 &    0.59 \\
  \hline
  \multirow{2}{*}[-0.3em]{\rotatebox[origin=c]{90}{$\text{med}_{50}$}} & inner &  3.84 &   3.58 &   3.22 &    2.44 &    3.55 &    3.23 &    2.86 &    2.04 &    2.96 &    2.63 &    2.10 &    1.14 \\ 
  & outer & 3.84 & 3.66 & 2.54 & 0.00 & 3.56 & 3.27 & 1.67 & 0.00 & 2.91 &    2.39 & 0.51 & 0.00 \\
\end{tabular}
\caption{Length of the estimated quantile-based $95\%$ credible intervals for the distributions of $M_{50}$ and $\text{med}_{50}$, induced by the inner and out spike and slab models for $\sigma$-stable hNRMIs, for different values of $\zeta$ and $\sigma$.}
\label{tab:only}
\end{table}

We complete the first part of the study by investigating the distributions of the number of observations in a sample coinciding with $x_0$, induced by the inner and the outer spike and slab models. The notation $N_0^{(m)}$, already introduced for the inner spike and slab model, will be henceforth used also to denote the same quantity induced by the outer spike and slab model. While for the outer model specification, it is easily verified that $N_0^{(m)}$ is a Binomial with parameters $m$ and $\zeta$, the distribution of $N_0^{(m)}$ induced by the inner model is provided in Theorem \ref{thm:N0_norm_ss} for the general case of hNRMIs, and coincides with \eqref{eq:N0_sigma_ss} when a $\sigma$-stable hNRMI spike and slab model is considered. Alternatively, the same distributions can be estimated by resorting to the samples generated in the first part of this section. We follow the latter approach and exclude from our analysis the case $\zeta=0$, given that it corresponds to $N_0^{(m)}\equiv0$. 

\begin{figure}[h!]
\centering
\includegraphics[width=0.33\textwidth,bb=250 10 990 710,clip=true]{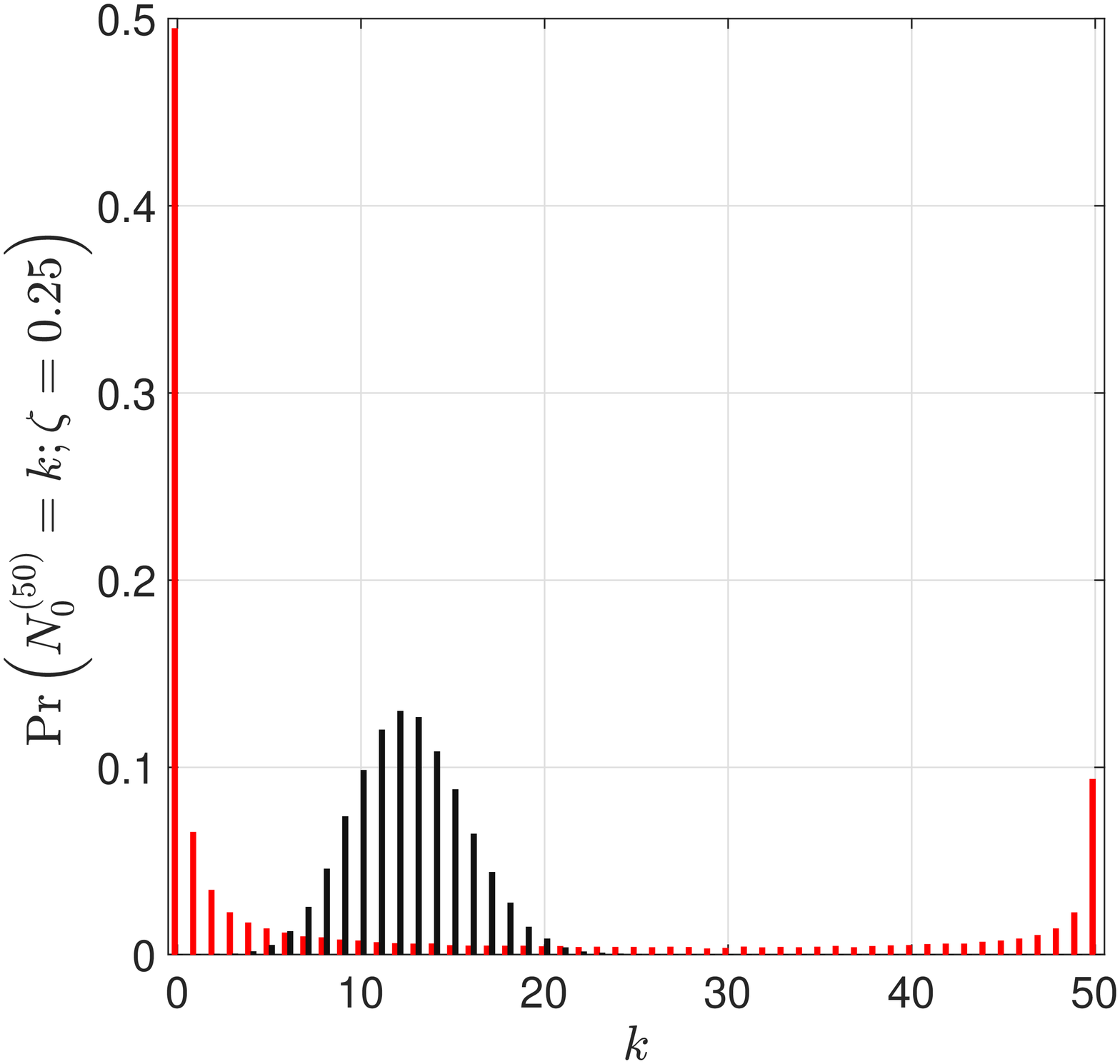}
\hspace{-0.1cm}%
\includegraphics[width=0.33\textwidth,bb=250 10 990 710,clip=true]{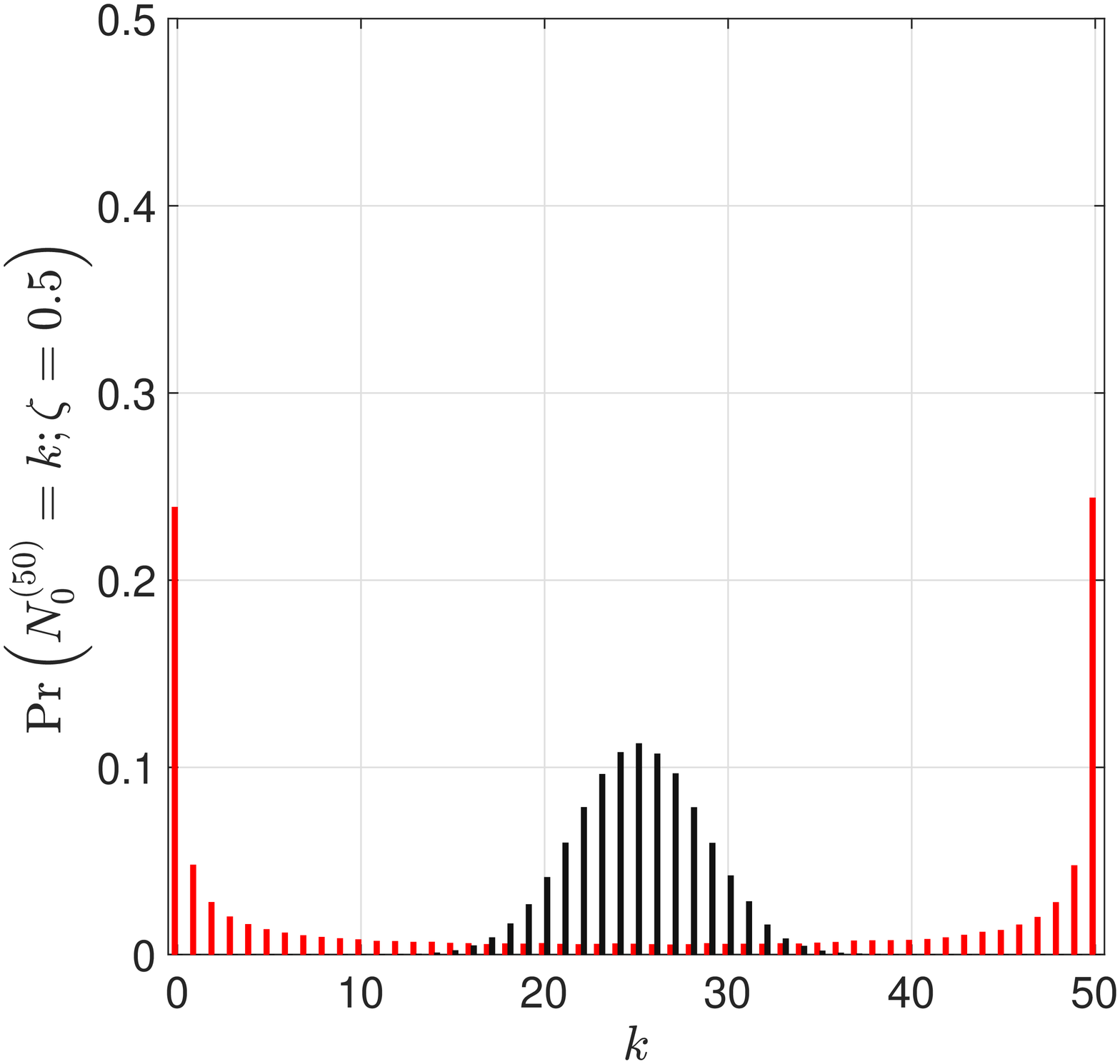}
\hspace{-0.1cm}%
\includegraphics[width=0.33\textwidth,bb=250 10 990 710,clip=true]{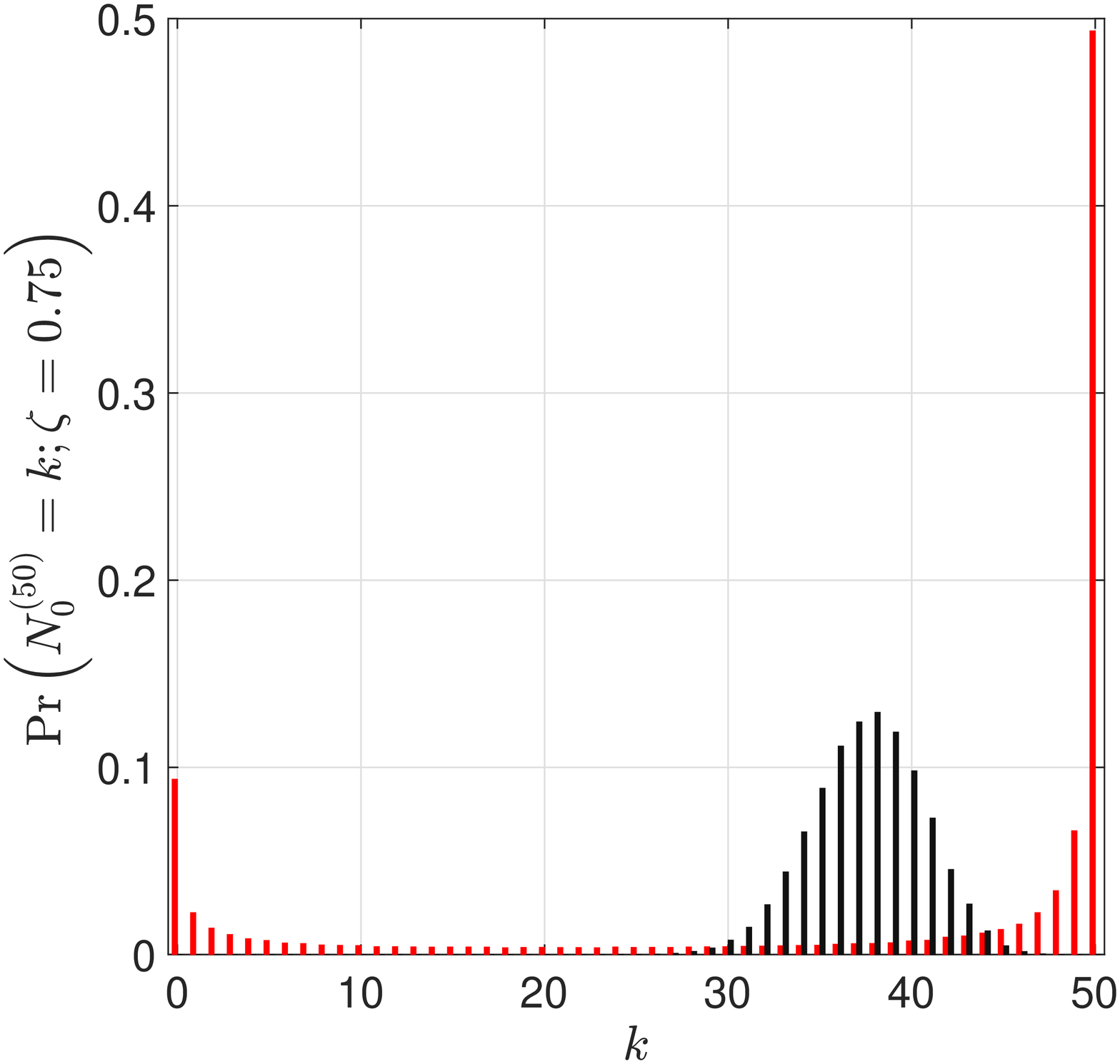}\\[5pt]
\includegraphics[width=0.33\textwidth,bb=250 10 990 710,clip=true]{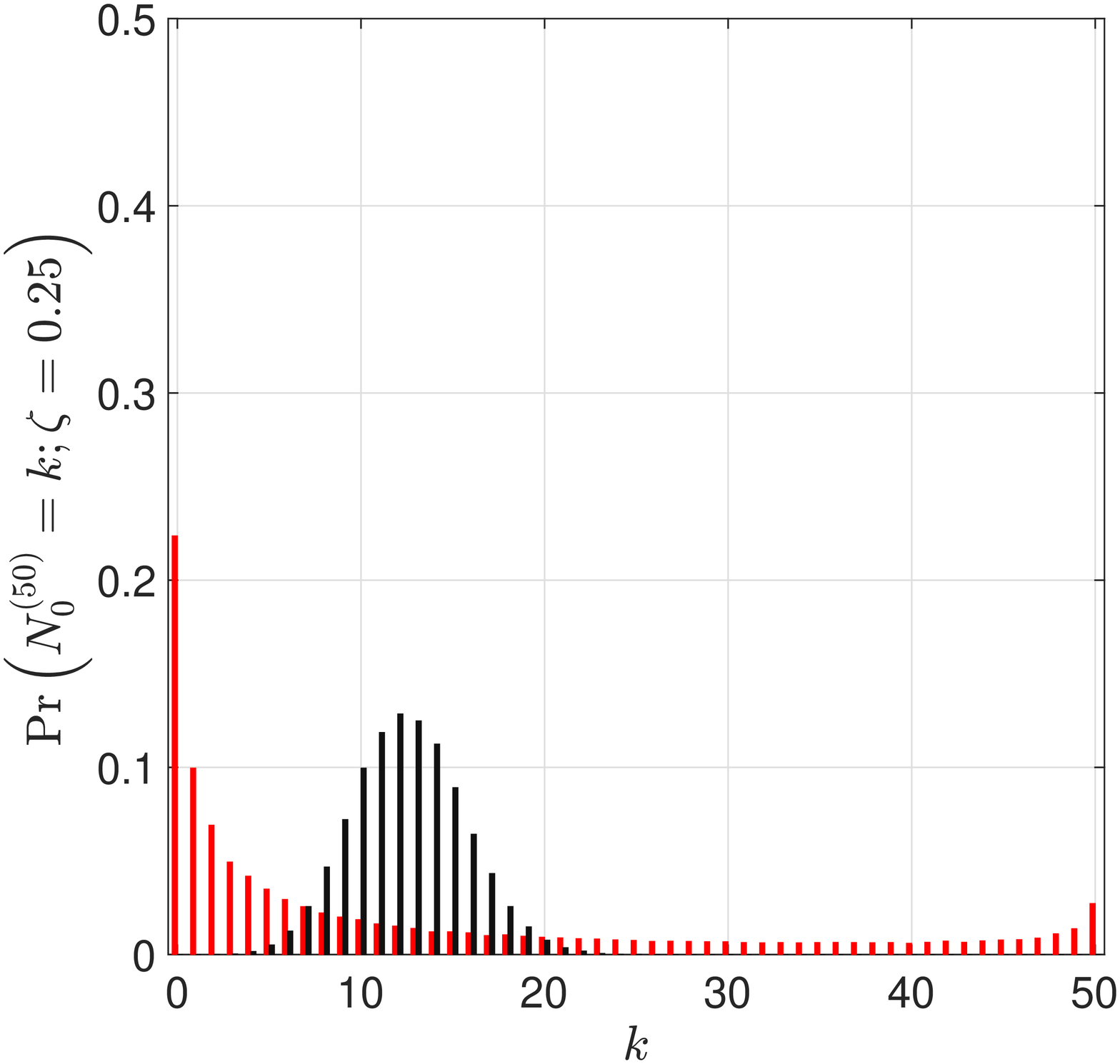}
\hspace{-0.1cm}%
\includegraphics[width=0.33\textwidth,bb=250 10 990 710,clip=true]{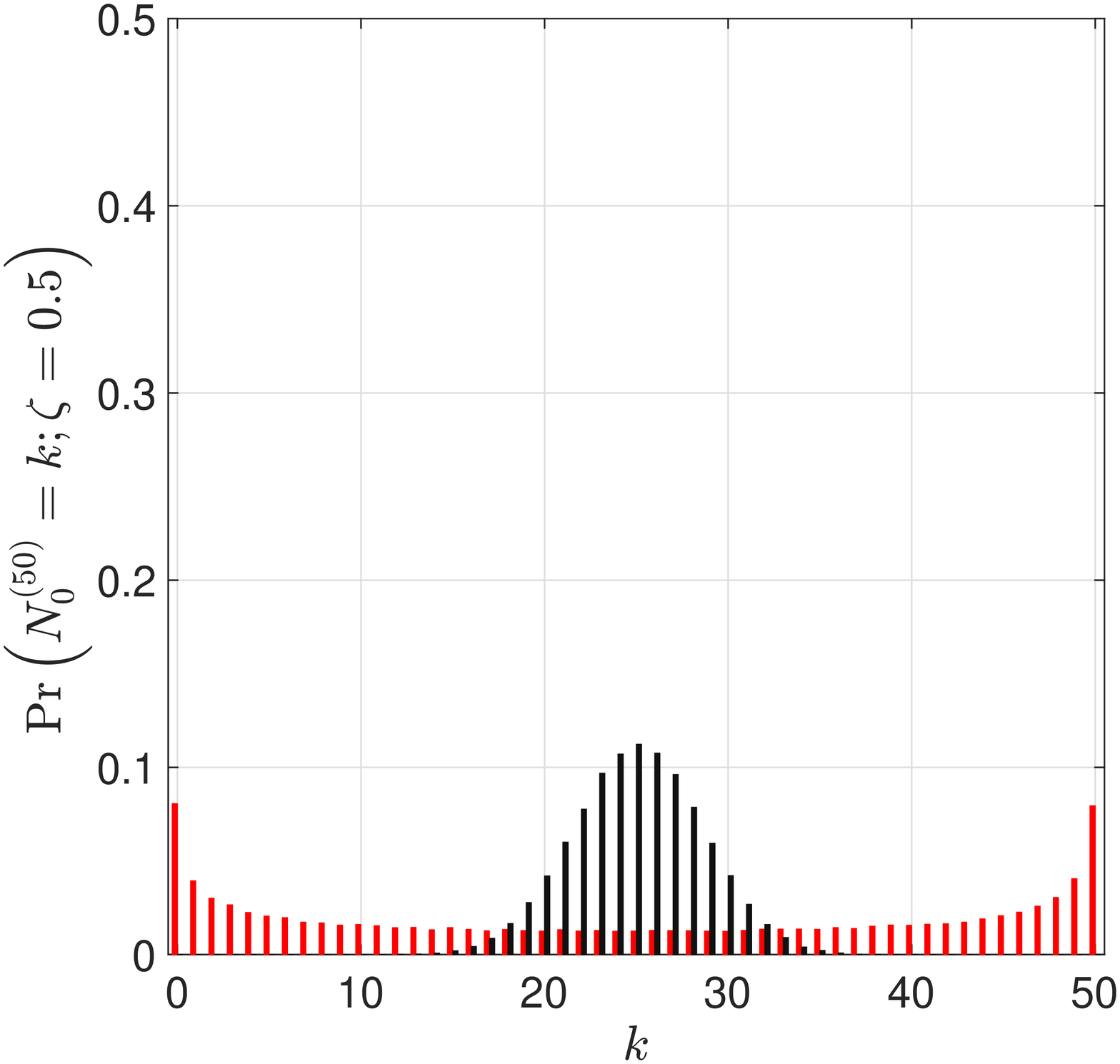}
\hspace{-0.1cm}%
\includegraphics[width=0.33\textwidth,bb=250 10 990 710,clip=true]{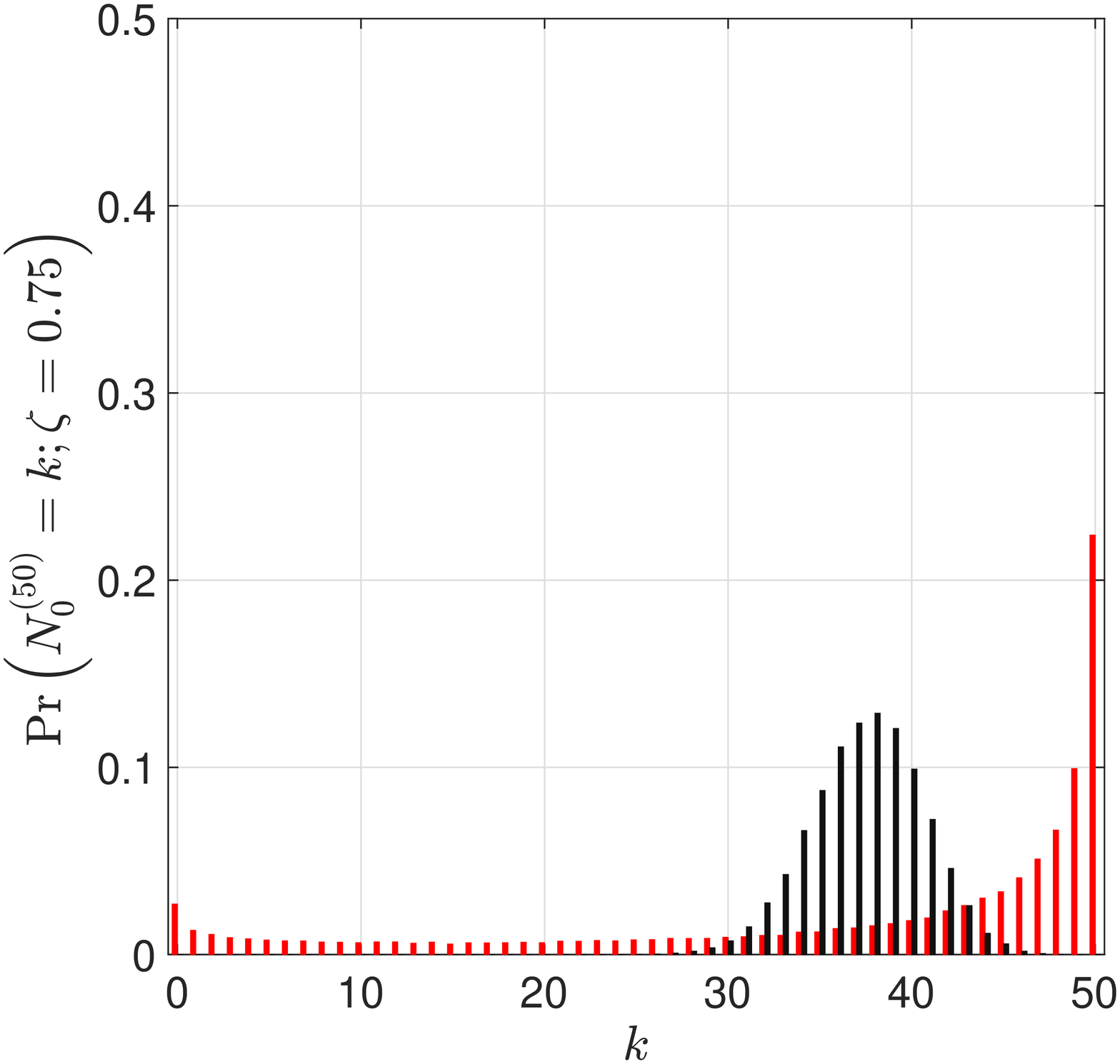}\\[5pt]
\includegraphics[width=0.33\textwidth,bb=250 10 990 710,clip=true]{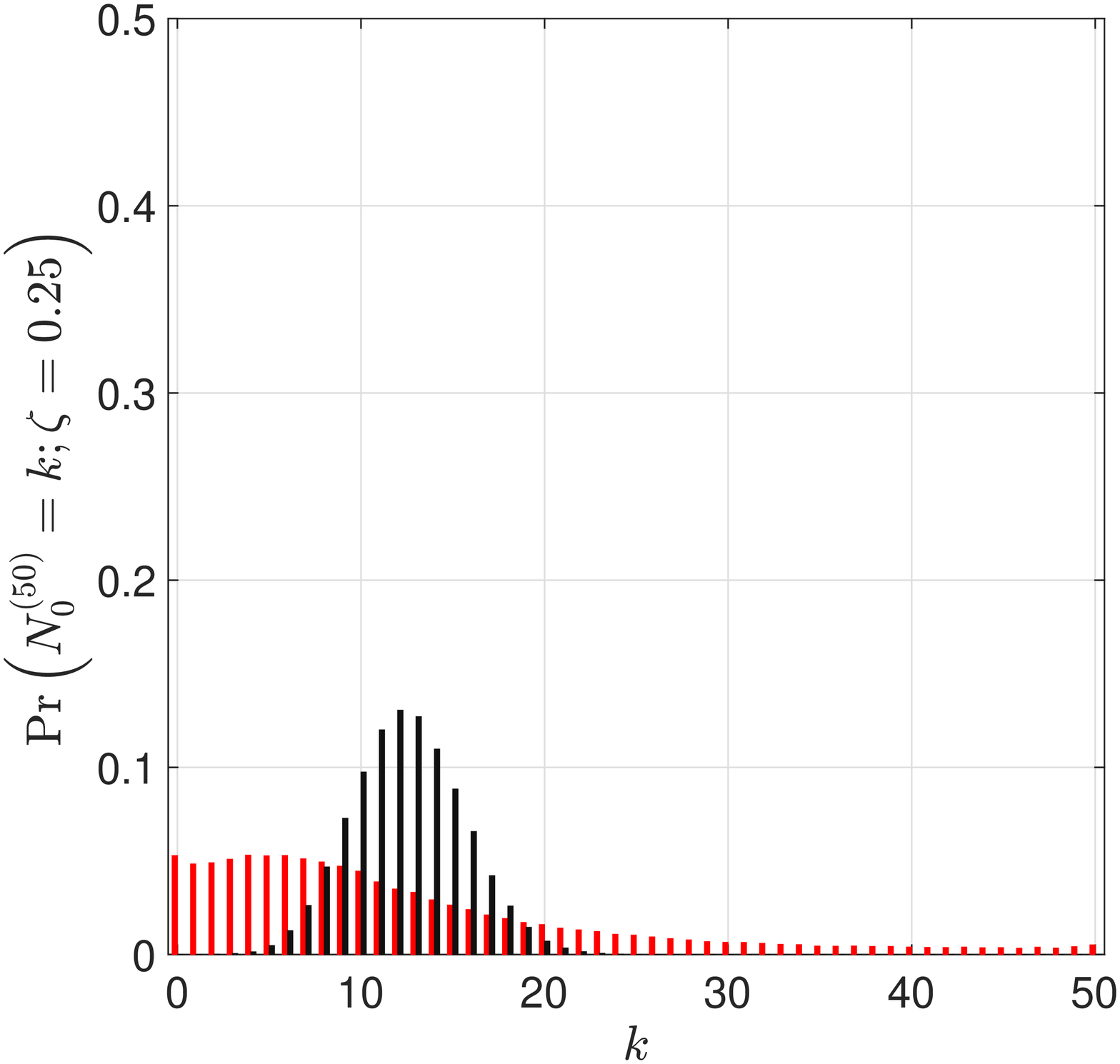}
\hspace{-0.1cm}%
\includegraphics[width=0.33\textwidth,bb=250 10 990 710,clip=true]{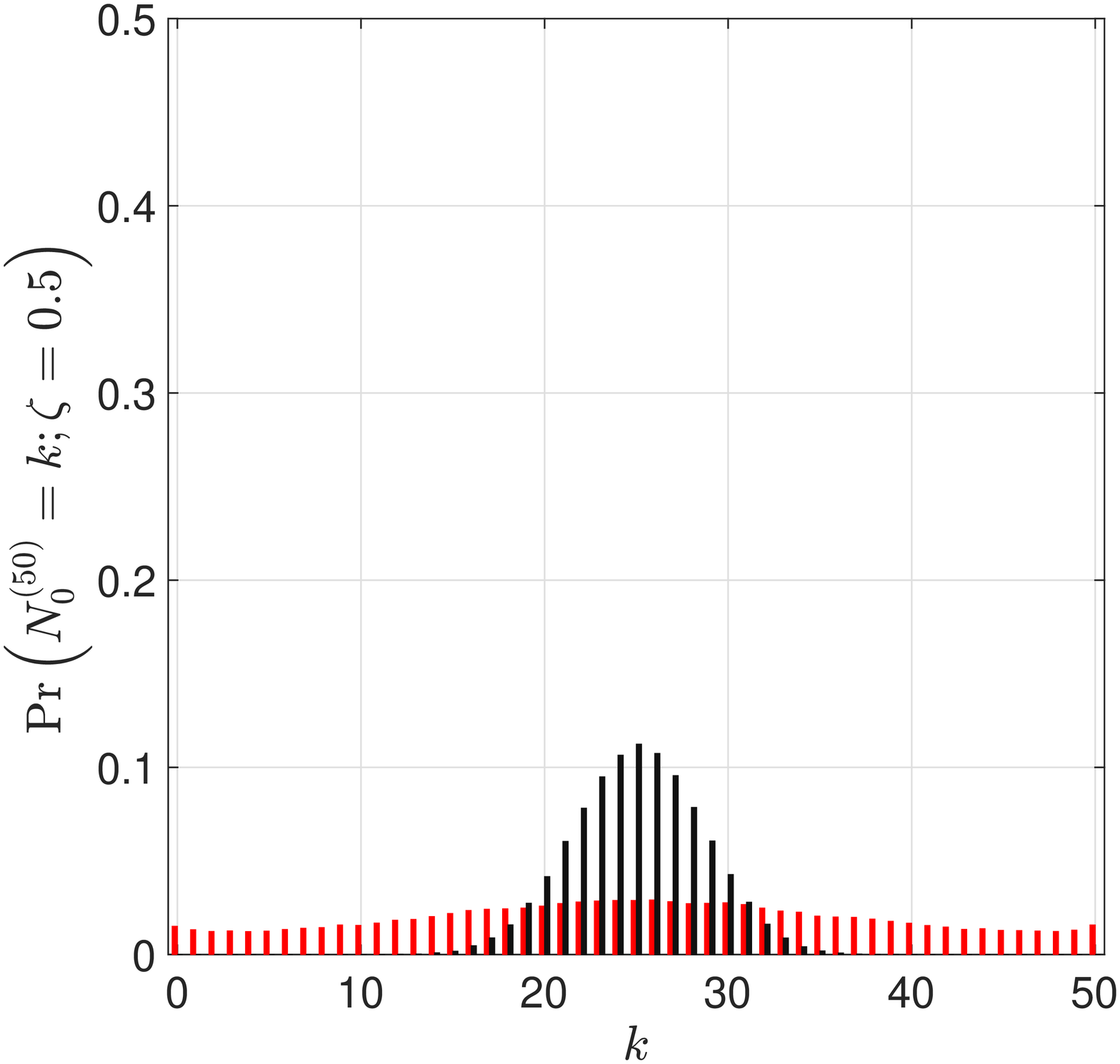}
\hspace{-0.1cm}%
\includegraphics[width=0.33\textwidth,bb=250 10 990 710,clip=true]{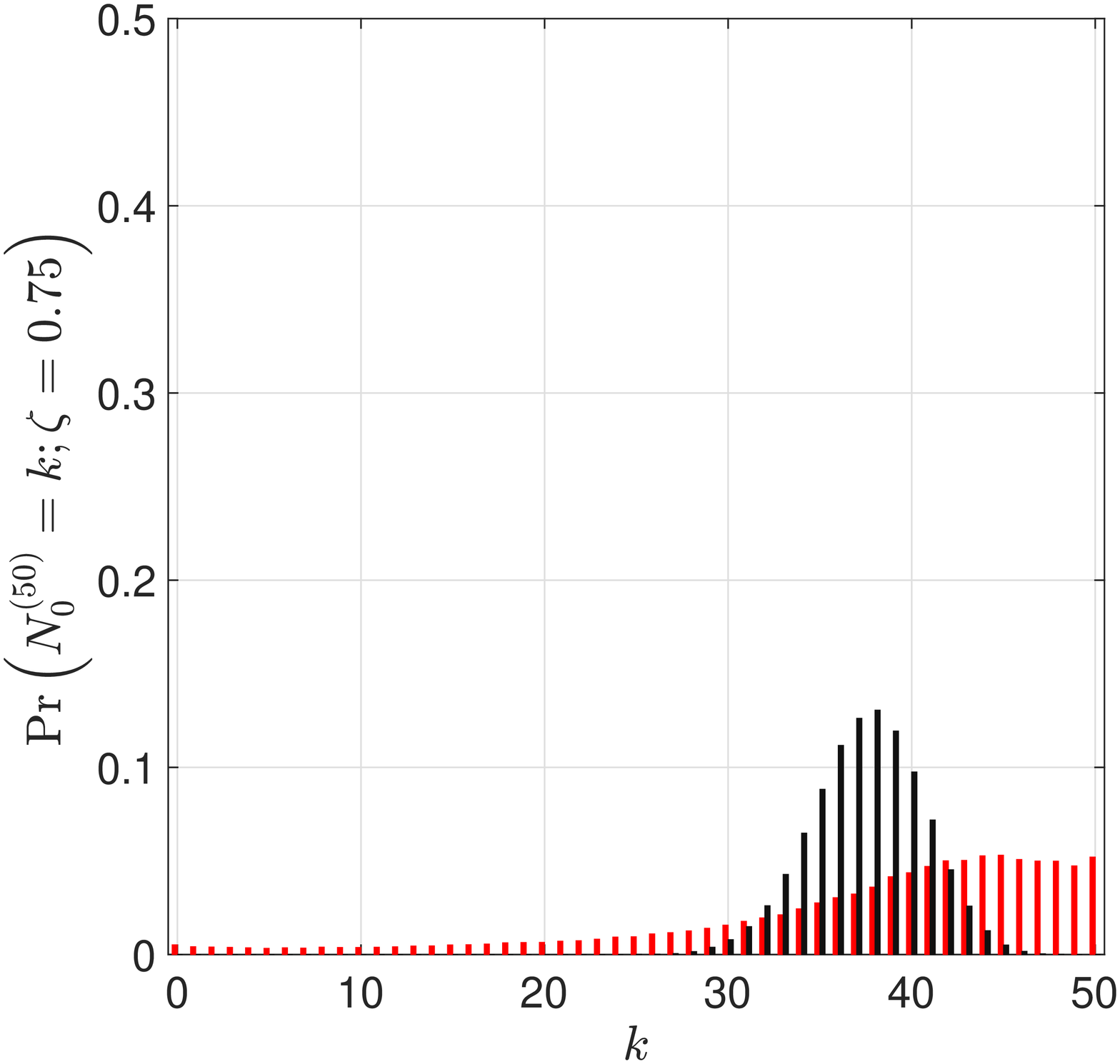}
\caption{Comparison of the distributions of $N_0^{(m)}$, for sample of size $m=50$, induced by inner (in red) and outer (in black) spike and slab models for $\sigma$-stable hNRMIs. Different rows refer to different values of $\sigma\in\{0.25,0.5,0.75\}$ (from top to bottom), different columns refer to different values of $\zeta\in\{0.25,0.5,0.75\}$ (from left to right).}\label{fig:N0}
\end{figure}

Figure \ref{fig:N0} displays a comparison between the distributions of $N_0^{(m)}$ induced by the two models, for $\zeta\in\{0.25,0.5,0.75\}$ and $\sigma\in\{0.25,0.5,0.75\}$, and with $m=50$. It is apparent, across all the combinations of the parameter values we considered, that the distribution of $N_0^{(m)}$ induced by $\tilde P$ is characterized by a larger variability than the one induced by $\tilde Q$. While for the latter most of the probability mass is allocated in a neighbourhood of the expected value, for the former the probability mass appears more spread over the whole support. Moreover, the study suggests that, for small values of $\sigma$, the inner spike and slab model, unlike its outer counterpart, assigns large probabilities to the extreme events $\{N_0^{(m)}=0\}$ and $\{N^{(m)}_0=m\}$. These distributional differences can be relevant when performing Bayesian classification by means of either spike and slab hNRMI prior model, and are in line with the findings of \citet{Can17} for the PY case.\\

Our numerical investigation of the distribution of $\tilde F_m$, its functionals $M_m$ and $\text{med}_m$, 
and the random variable $N_0^{(m)}$, although limited to the case of $\sigma$-stable hNRMIs, underpins the interpretation we gave of Proposition \ref{prop:var}. 
Our study indicates that, if compared to the inner model, the outer one assigns less prior probability to realizations of $\tilde F_{m}$ that deviate considerably from $F_0$, as well as to values of $N_0^{(m)}$ which deviate considerably from its expected value. The inner spike and slab model appears characterized by a larger variability and thus can be interpreted as less informative than the outer one. In most contexts this represents a desirable feature and encapsulates the advantage of using a fully nonparametric specification as the inner spike and slab model is.

\subsection{Inner vs outer models a posteriori}

\begin{figure}[h!]
\centering
\includegraphics[width=0.33\textwidth,bb=230 10 990 710,clip=true]{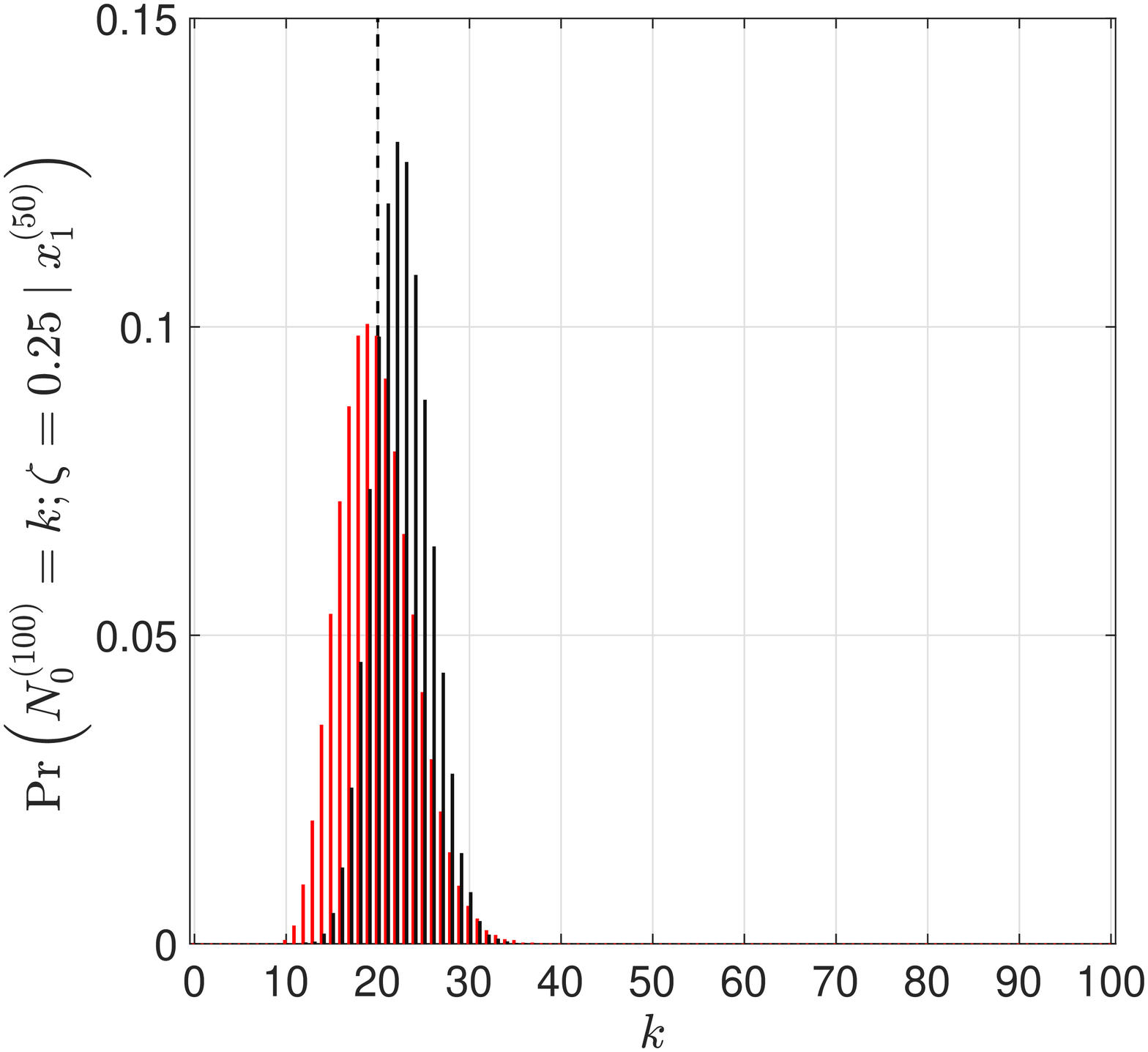}
\hspace{-0.1cm}%
\includegraphics[width=0.33\textwidth,bb=230 10 990 710,clip=true]{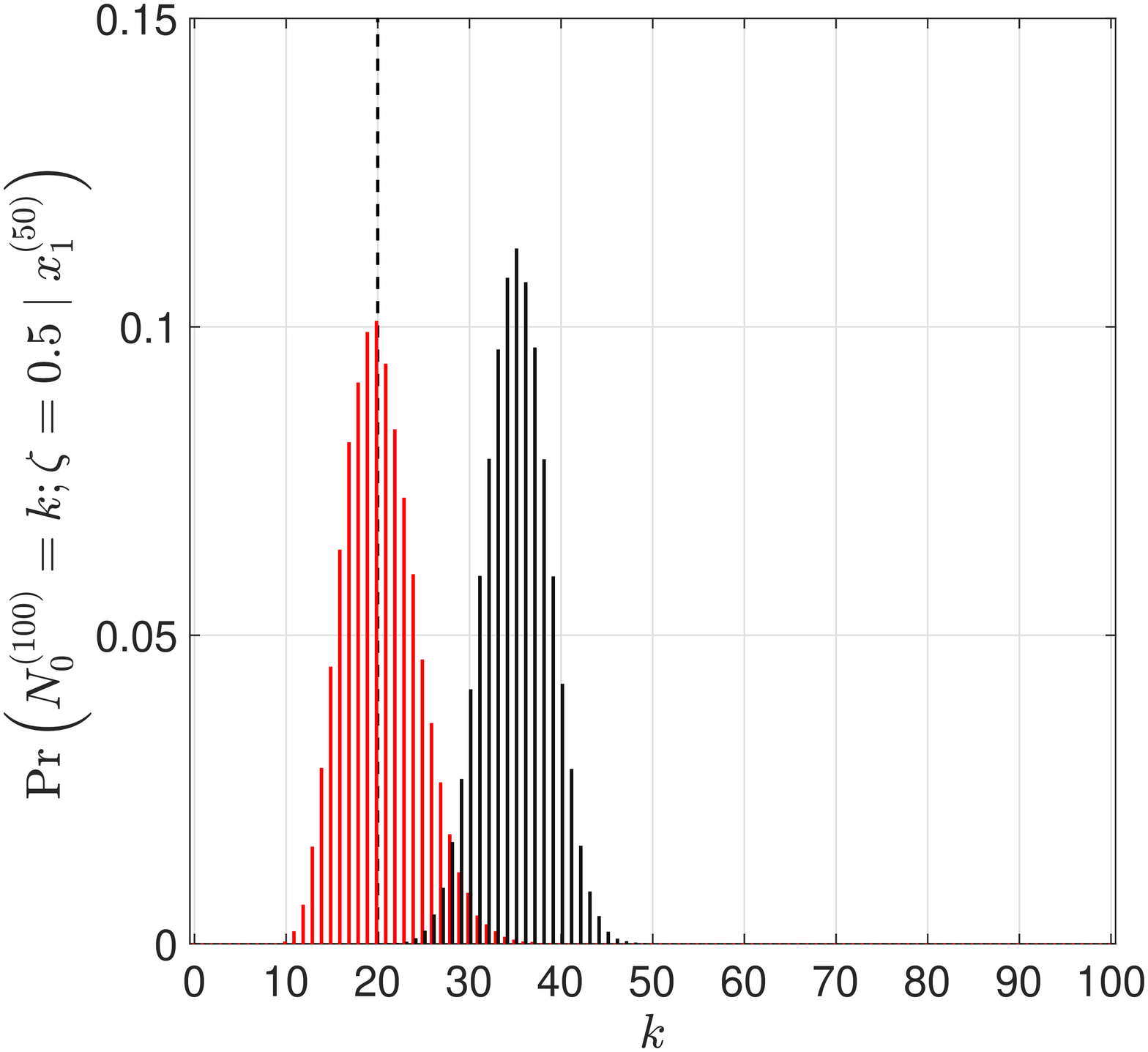}
\hspace{-0.1cm}%
\includegraphics[width=0.33\textwidth,bb=230 10 990 710,clip=true]{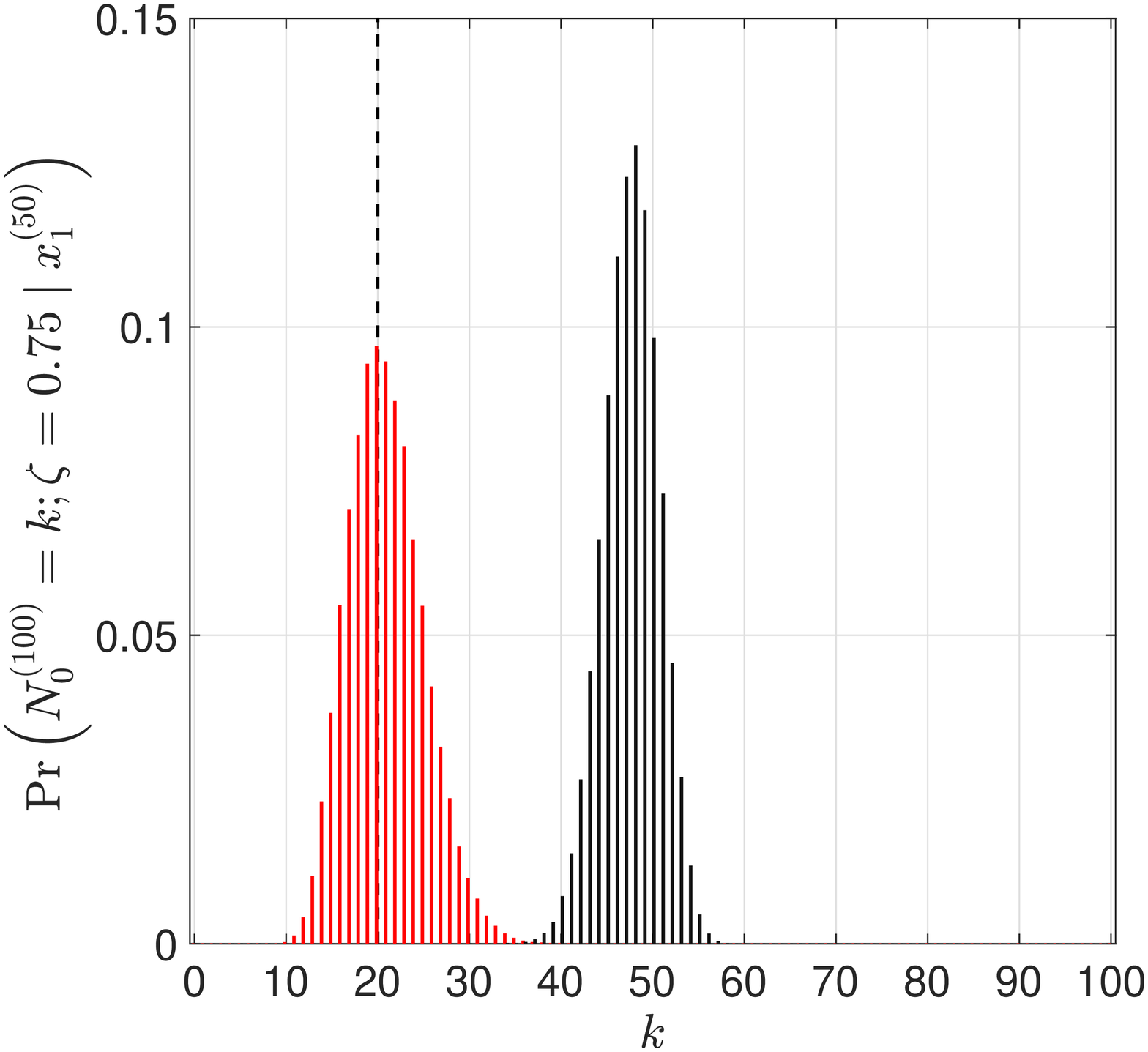}\\[5pt]
\includegraphics[width=0.33\textwidth,bb=230 10 990 710,clip=true]{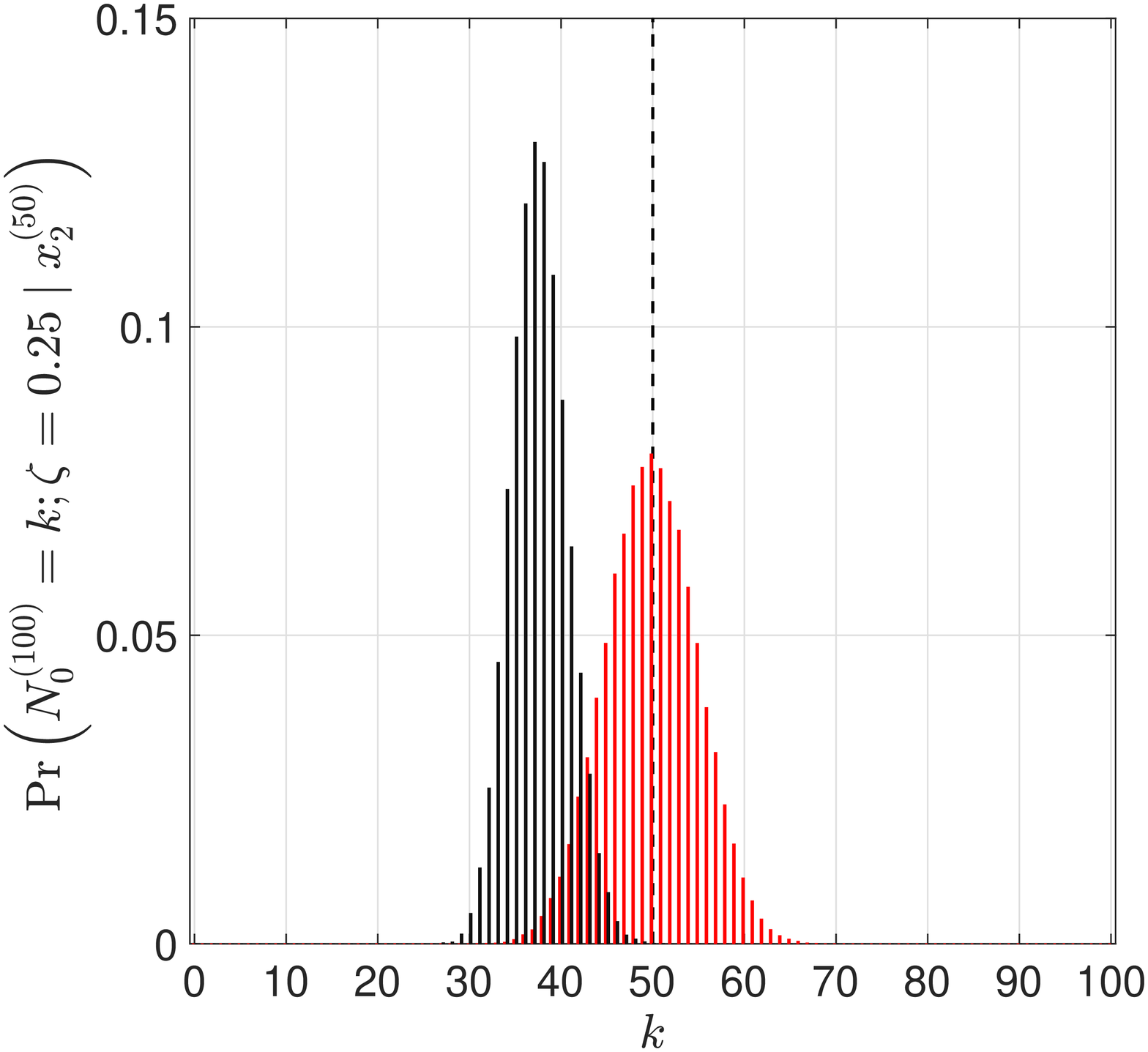}
\hspace{-0.1cm}%
\includegraphics[width=0.33\textwidth,bb=230 10 990 710,clip=true]{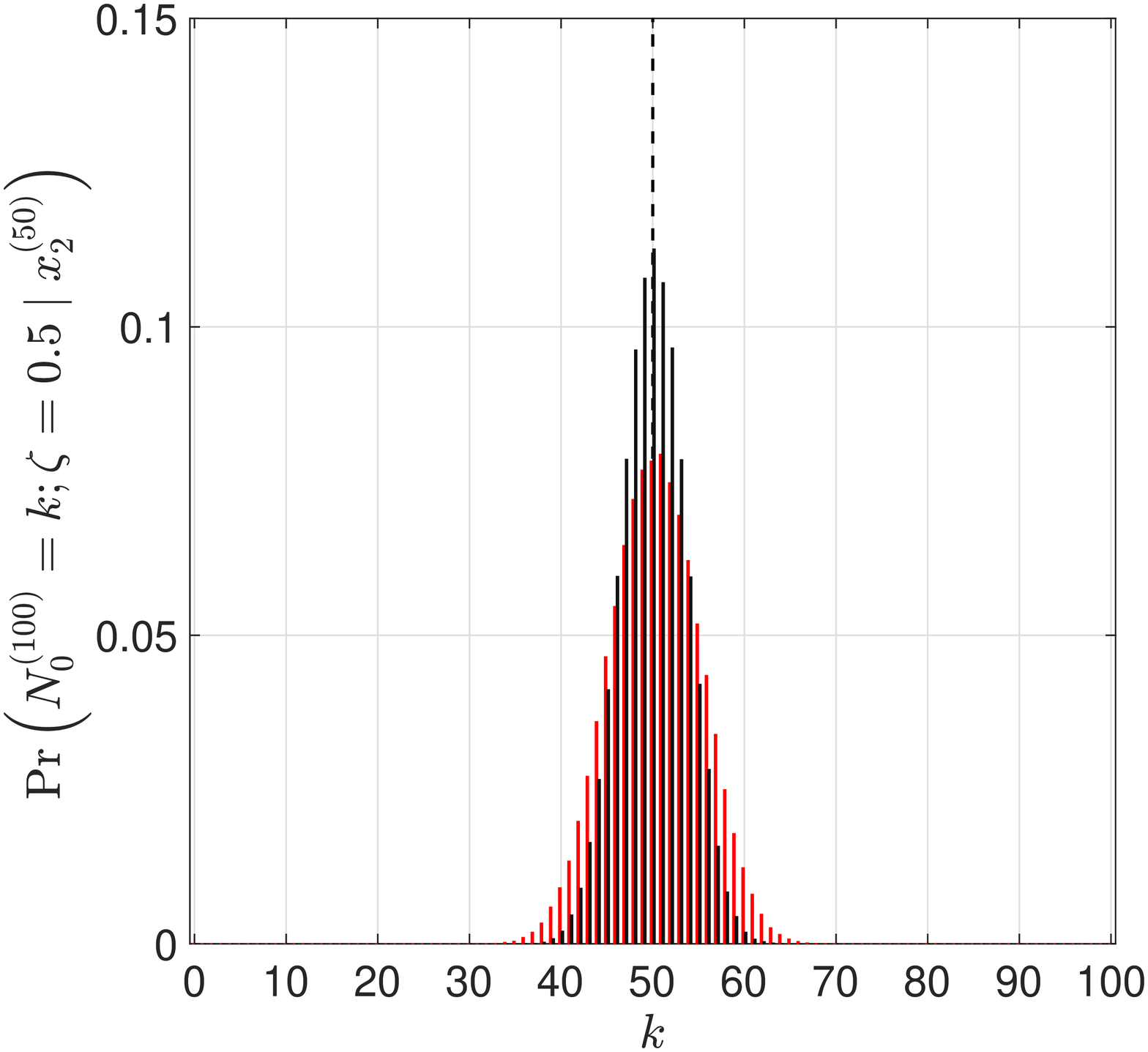}
\hspace{-0.1cm}%
\includegraphics[width=0.33\textwidth,bb=230 10 990 710,clip=true]{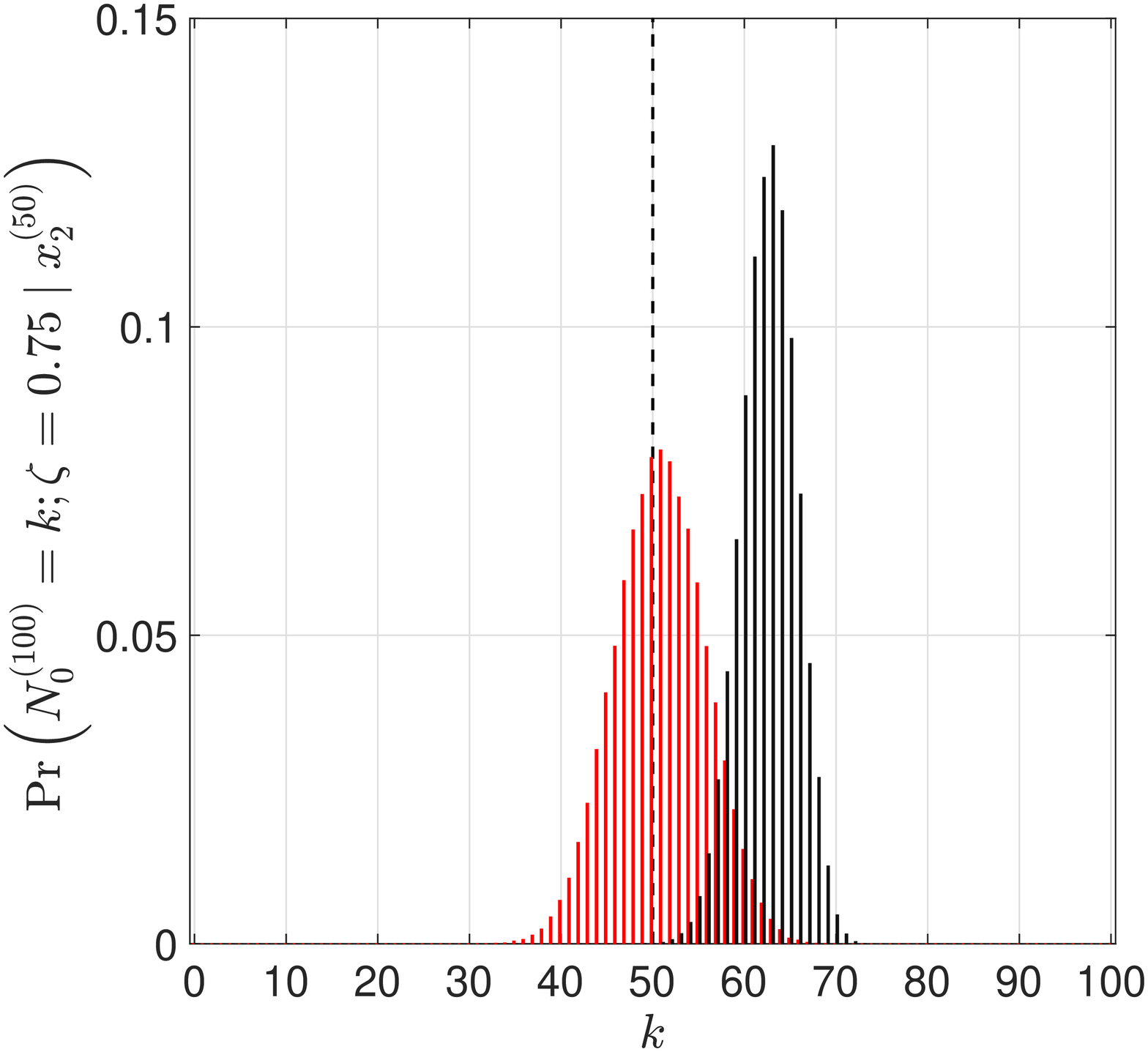}\\[5pt]
\includegraphics[width=0.33\textwidth,bb=230 10 990 710,clip=true]{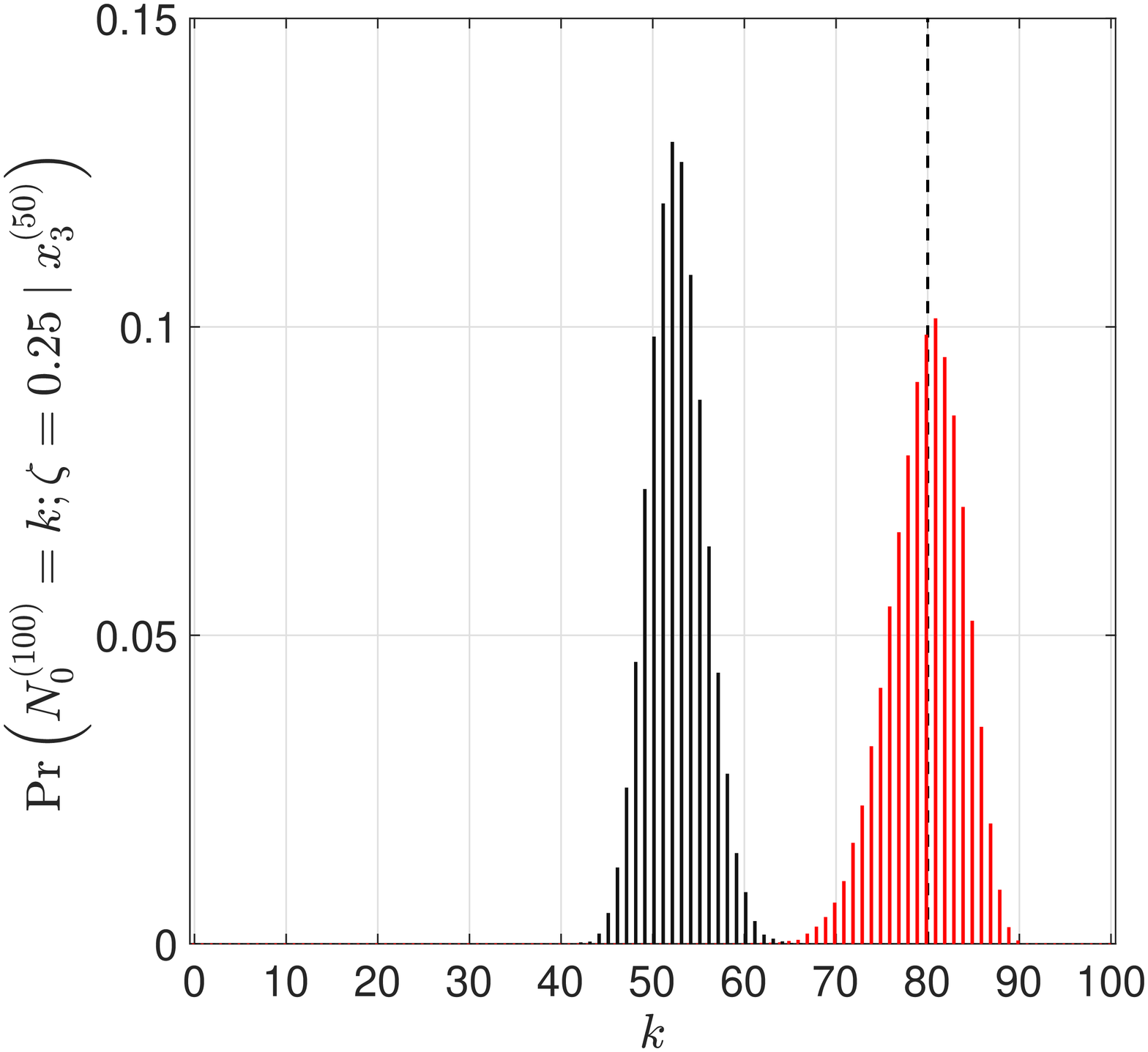}
\hspace{-0.1cm}%
\includegraphics[width=0.33\textwidth,bb=230 10 990 710,clip=true]{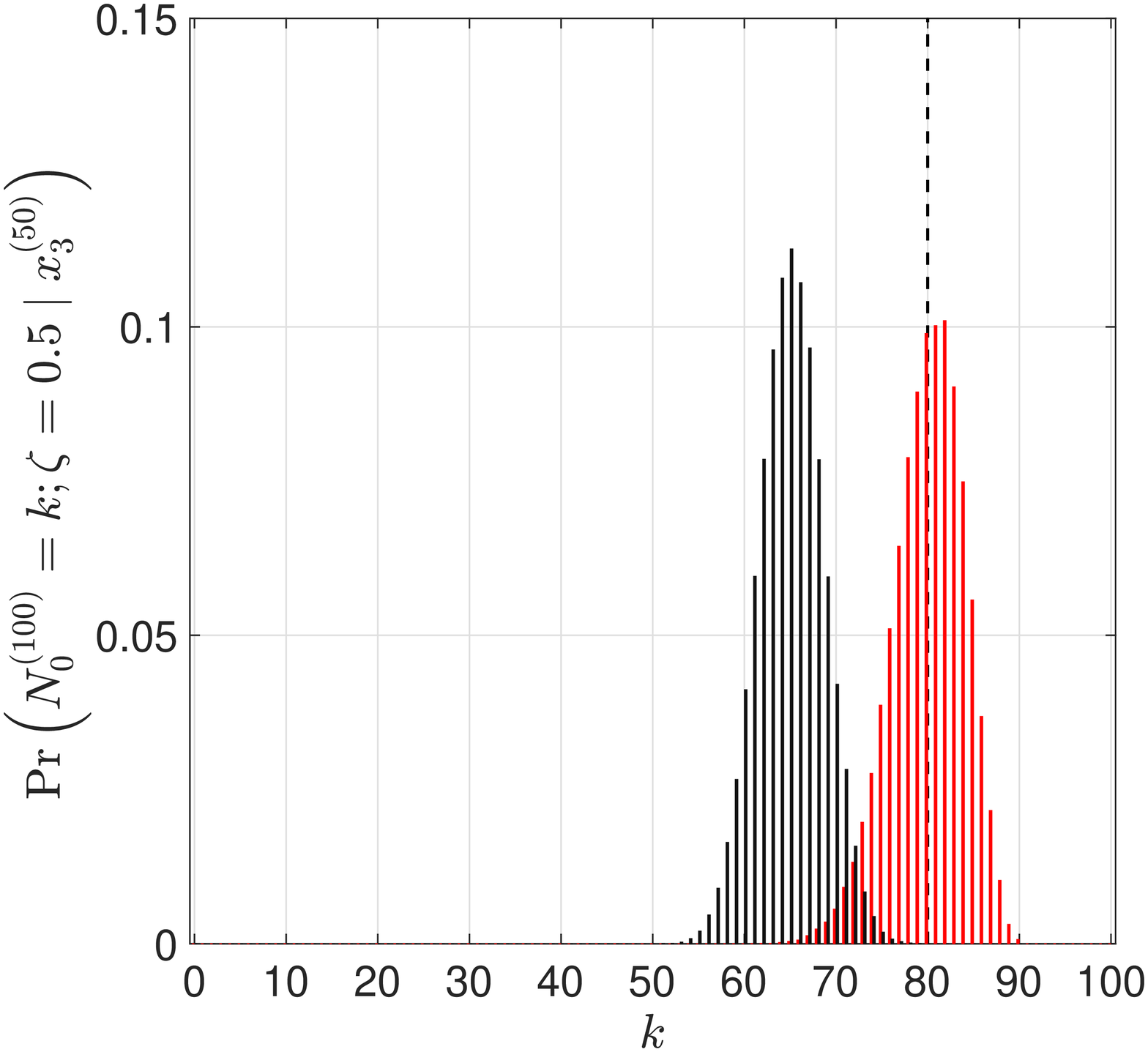}
\hspace{-0.1cm}%
\includegraphics[width=0.33\textwidth,bb=230 10 990 710,clip=true]{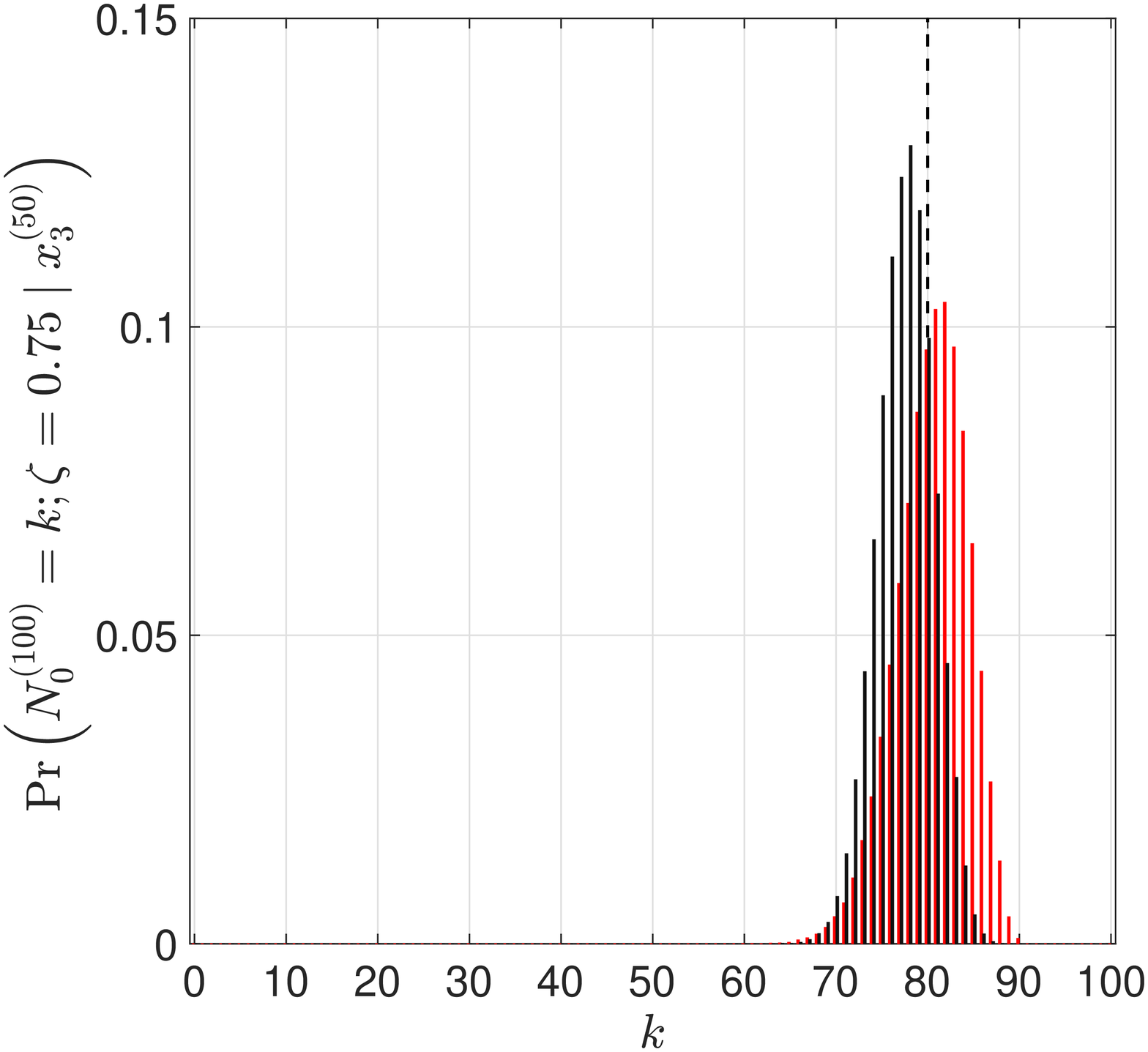}
\caption{Comparison of the distributions of $N_0^{(n+m)}$, with $n=m=50$, induced by inner (in red) and outer (in black) spike and slab models for $\sigma$-stable hNRMIs, with $\sigma=0.25$. Different rows refer to different samples $x_\ell^{(50)}$, with $\ell\in\{1,2,3\}$ (from top to bottom), different columns refer to different values of $\zeta\in\{0.25,0.5,0.75\}$ (from left to right). The vertical dashed black line indicates, in each panel, the value $\pi_\ell (n+m)$, that is the number of observations in the enlarged sample corresponding to the sample proportion $\pi_\ell$.}\label{fig:N0_post}
\end{figure}

We next investigate the posterior behaviour of inner and outer spike and slab models, with a focus on the posterior distribution of $N_0^{(n+m)}$, conditional on the observation of a sample $x^{(n)}$. Such investigation allows us to better understand how the prior specification for the two models affects the posterior probability of observations coinciding with $x_0$. As in Section \ref{sec:priorstudy}, we focus on the case of $\sigma$-stable hNRMI models. We set $n=50$ and consider three observed samples, denoted as $x_1^{(50)}$, $x_2^{(50)}$ and $x_3^{(50)}$. 
In $x_{1}^{(50)}$, $10$ observations coincide with $x_0$ (sample proportion $\pi_1=0.2$) and the remaining 40 display 7 distinct values with frequencies $(25,5,3,3,2,1,1)$; in
   $x_2^{(50)}$, $25$ observations coincide with $x_0$ (sample proportion $\pi_2=0.5$) and the remaining 25 display 7 distinct values with frequencies $(10,5,3,3,2,1,1)$;  
 in $x_3^{(50)}$, $40$ observations coincide with $x_0$ (sample proportion $\pi_3=0.8$) and the remaining 10 display 5 distinct values with frequencies $(3,3,2,1,1)$. 
Conditionally on each of the three samples, we generate realizations $(x_{n+1},\ldots,x_{n+m})$ of $(X_{n+1},\ldots,X_{n+m})$ by considering $\sigma\in\{0.25,0.5,0.75\}$ and $\zeta\in\{0.25,0.5, 0.75\}$. Specifically, we set $m=50$ and, for each $x_\ell^{(n)}$, with $\ell=1,2,3$, and each combination of values for $\sigma$ and $\zeta$, we generate $100\,000$ samples and use them to estimate the posterior distribution of $N_0^{(n+m)}$, conditional on $x_\ell^{(n)}$. As for the outer spike and slab model, samples are generated by using Algorithm \ref{alg:PUS}, as done in Section \ref{sec:study1}. Figure \ref{fig:N0_post} refers to the case $\sigma=0.25$ and displays a structural difference between inner and outer spike and slab models. Although not displayed here, similar conclusions can be drawn by investigating the cases $\sigma=0.5$ and $\sigma=0.75$. The posterior distribution of $N_0^{(n+m)}$, conditional on the sample $x_\ell^{(n)}$, for $\ell=1,2,3$, appears concentrated around $\pi_\ell(n+m)$ when induced by the inner spike and slab model, regardless of the value of $\zeta$. The same does not happen for the posterior distribution of $N_0^{(n+m)}$ induced by the outer spike and slab model, which is instead highly sensitive to the value of $\zeta$. In other terms, as reported in Table \ref{tab:N0_post}, the expected proportion of observations coinciding with $x_0$ in the enlarged sample of size $n+m$,  is close to the observed proportion $\pi_\ell$ when an inner spike and slab model is adopted; on the contrary, the outer spike and slab model leads to an expected proportion of observations coinciding with $x_0$ in $X^{(n+m)}$ which is shrunk towards the expected prior proportion $\zeta$.
Table \ref{tab:N0_post} shows that this behaviour is common across all the considered values of $\sigma$. The empirical results of our study indicates that the prior induced on $N_0^{(m)}$ by the inner spike and slab model is less informative than the one induced by the outer spike and slab model. Our findings also indicate that, while a hyperprior for $\zeta$ seems needed when the outer spike and slab model is adopted, the inner specification appears structurally robust to any specification of $\zeta$.

\begin{table}[h!]
\makegapedcells
\centering
\begin{tabular}{ c c | c c c | c c c | c c c }
\multicolumn{2}{c|}{$\sigma$} & \multicolumn{3}{c|}{$0.25$} & \multicolumn{3}{c|}{$0.5$} & \multicolumn{3}{c}{$0.75$} \\
  \hline
  \multicolumn{2}{c|}{$\zeta$} &  $0.25$ &  $0.5$ &  $0.75$ &  $0.25$ &  $0.5$ &  $0.75$ &  $0.25$ &  $0.5$ &  $0.75$\\
  \hline
  \multirow{2}{*}[-0.3em]{\rotatebox[origin=c]{0}{$\pi_1=0.2$}} & inner & 0.20 & 0.20 & 0.21 & 0.20 &  0.21 & 0.22 & 0.22 & 0.20 & 0.22 \\
  & outer & 0.23 & 0.35 & 0.48 & 0.22 & 0.35 &    0.48 & 0.22 & 0.35 & 0.47 \\
    \hline
  \multirow{2}{*}[-0.3em]{\rotatebox[origin=c]{0}{$\pi_2=0.5$}} & inner &    0.50 &   0.50 & 0.51 & 0.50 & 0.50 &  0.51 & 0.51 &  0.49 &  0.51  \\
  &outer & 0.38 & 0.50 &  0.63 & 0.37 & 0.50 &   0.62 & 0.37 & 0.50 & 0.63 \\
    \hline 
  \multirow{2}{*}[-0.3em]{\rotatebox[origin=c]{0}{$\pi_3=0.8$}} & inner &    0.80 &  0.80 & 0.81 & 0.80 & 0.80 &  0.80 &  0.80 & 0.79 & 0.80  \\
  &outer & 0.53 & 0.65 & 0.77 &  0.52 &  0.65 &   0.78 & 0.52 & 0.65 & 0.78 \\
\end{tabular}
\caption{
	Estimated proportion of observations coinciding with $x_0$ in the sample $X^{(100)}$, conditional on the observation of $x_1^{(50)}$ ($\pi_1=0.2$), $x_2^{(50)}$ ($\pi_2=0.5$), and $x_3^{(50)}$ ($\pi_3=0.8$), for values of $\sigma\in\{0.25,0.5,0.75\}$ and $\zeta\in\{0.25,0.5,0.75\}$.}
\label{tab:N0_post}
\end{table}

\section{Proofs}\label{sec:proofs}
This section contains the proofs of the main results of Sections \ref{sec:nrmi} and \ref{sec:main}.
\subsection{Proof of Proposition \ref{prop:var}}
The proof follows from a simple application of Theorem \ref{eppf:nrmi}. Indeed, let
$p=\Pi_1^{(2)}(2;0)=1-\Pi_{2}^{(2)}(1,1;0)$
and note that 
\[
\E \left[\tilde P^2(f)\right]=\E\left[ \int_\X f^2(x)\,\tilde P^2(\ddr x)\right]+\E\left[\int_{\X^2_*}f(x_1)f(x_2)\,\tilde P(\ddr x_1)\,\tilde P(\ddr x_2)\right],
\]
where $\X^2_*=\{(x_1,x_2)\in\X^2:\: x_1\ne x_2\}$. Since $\xi_{2,2}(u)=\tau_2(u)$ and $\xi_{2,1}(u)= \tau_1^2(u)$, from \eqref{eq:eppf:nrmi} one has
\begin{align*}
\E\left[\int_\X f^2(x)\,\tilde P^2(\ddr x)\right]
&=p(1-\zeta)P^*(f^2)+p\zeta f^2(x_0)+(1-p)\zeta^2f^2(x_0),\\[4pt]
\E\left[\int_{\X^2_*}f(x_1)f(x_2)\,\tilde P(\ddr x_1)\,\tilde P(\ddr x_2)\right]
&= (1-p)(1-\zeta)^2\,\Big(P^*(f)\Big)^2+2(1-p)\zeta(1-\zeta)f(x_0)P^*(f).
\end{align*}
By combining these two, one obtains
\[
\mbox{var}(\tilde P(f))=
p(1-\zeta)P^*(f^2)+p\zeta(1-\zeta)f^2(x_0)-p(1-\zeta)^2\left(P^*(f)\right)^2-2p\zeta(1-\zeta)f(x_0)P^*(f).
\]
Moreover, using \eqref{eq:eppf_hnrmi} it can be easily seen that 
\[
\mbox{var}(\tilde Q(f))=(1-\zeta)^2\,\left\{p\,P^*(f^2)+(1-p)\,\left(P^*(f)\right)^2 \right\}-(1-\zeta)^2\left(P^*(f)\right)^2
\]
and the result follows. \hfill \qed
\subsection{Proof of Theorem \ref{eppf:nrmi}}

The proof relies on techniques similar to those used in the proof of Theorem 2 in \citet{Can17}. By definition, we have
\begin{equation}\label{eq:EPPF_int}
    \Pi_k^{(n)}(n_1,n_2,\ldots,n_k;\zeta)=\E\left[\int_{\X^{k}}\prod_{j=1}^{k}\tilde P(\ddr x_j)\right],
\end{equation}
for any vector of positive integers $(n_1,n_2,\ldots,n_k)$ such that $\sum_{j=1}^k n_k=n$, and where the integrating variables are such that $x_1\neq x_2 \neq \ldots \neq x_k$. The right-hand side of \eqref{eq:EPPF_int} can be written as 
\begin{equation}\label{eq:EPPF_gamma_int}
    \frac{1}{\Gamma(n)}\int_{\X^k} \int_0^\infty u^{n-1} \E\left[ \edr^{-u\tilde\mu(\X)}\prod_{j=1}^k \tilde\mu^{n_j}(\ddr x_j)\right]\ddr u,
\end{equation}
where we used Fubini's theorem.  Henceforth we let $\d x_j$ denote, for any $j=1,\ldots,k$, a neighbourhood of radius $\varepsilon$ of $x_j$ for which, thus, $\d x_j \downarrow \{x_j\}$ as $\varepsilon\downarrow 0$. We then start by focusing on the expected value appearing in \eqref{eq:EPPF_gamma_int}, for which we consider two cases: i) none of the infinitesimal intervals $\ddr x_j$ contains $x_0$, ii) $x_0\in \ddr x_l$ for one $l\in\{1,\ldots,k\}$. Note that, given the infinitesimal nature of the intervals $\ddr x_j$, with $j=1,2,\ldots,k$, and the fact that $x_1\neq x_2\neq\ldots\neq x_k$, $x_0$ is not contained in more than one interval $\ddr x_j$ with probability one. 
\begin{itemize}
    \item[i)] In this case, only the diffuse component $P^*$ in $P_0$ plays a role and it can be showed by standard techniques that
    \begin{align}\label{eq:first_case}
        \E\left[ \edr^{-u\tilde\mu(\X)}\prod_{j=1}^k \tilde\mu^{n_j}(\ddr x_j)\right]&\cong c^k(1-\zeta)^k \edr^{-c\psi(u)} \prod_{j=1}^k \tau_{n_j}(u) P^*(\ddr x_j),
    \end{align}
    where \eqref{eq:first_case} provides a first-order approximation of the expected value on the left-hand side, which is everything we need as the higher order terms vanish when computing the integral over $\X^k$ in \eqref{eq:EPPF_gamma_int}.
    \item[ii)] We define $\X^*=\X\setminus\{x_1,\ldots,x_k\}$ and exploit the independence of increments of $\tilde \mu$ to write 
    \begin{equation}\label{eq:prod_exp}
        \E\left[ \edr^{-u\tilde\mu(\X)}\prod_{j=1}^k \tilde\mu^{n_j}(\ddr x_j)\right]=\E\left[\edr^{-u\tilde \mu(\X^*)}\right]\prod_{m\neq l}\E\left[\e^{-u\tilde \mu(\ddr x_m)}\tilde \mu^{n_m}(\ddr x_m)\right]\E\left[\e^{-u\tilde \mu(\{x_0\})}\tilde \mu^{n_l}(\{x_0\})\right].
    \end{equation}
    Again, in the first two terms of the right-hand side of \eqref{eq:prod_exp}, only the diffuse component $P^*$ of $P_0$ contributes to the integral. Thus, it can be showed by standard techniques that
    \begin{align}\label{eq:first}
      \E\left[\edr^{-u\tilde \mu(\X^*)}\right]&=\edr^{-c(1-\zeta)P^*(\X^*)\psi(u)}
      \end{align}
      and
      \begin{align}\label{eq:second}
      \prod_{m\neq l}\E\left[\e^{-u\tilde \mu(\ddr x_m)}\tilde \mu^{n_m}(\ddr x_m)\right]&\cong c^{k-1}(1-\zeta)^{k-1}\prod_{m\neq l}\edr^{-c(1-\zeta)P^*(\ddr x_m)\psi(u)}P^*(\ddr x_m)\tau_{n_m}(u),
    \end{align}
    where the approximation can be interpreted as the one in \eqref{eq:first_case}
    As for the last term of \eqref{eq:prod_exp}, we apply Feynman's technique for integration and write
    \begin{align*}
       \E\left[\e^{-u\tilde \mu(\{x_0\})}\tilde \mu^{n_l}(\{x_0\})\right]&=(-1)^{n_l}\frac{\ddr^{n_l}}{\ddr u^{n_l}} \E\left[\edr^{-u\tilde \mu(\{x_0\})}\right]\\
       &=(-1)^{n_l} \frac{\ddr^{n_l}}{\ddr u^{n_l}} \edr^{-c\zeta\psi(u)}\\
       &=\sum_{i=1}^{n_l}c^i\zeta^i \edr^{-c\zeta \psi(u)}\frac{1}{i!}\sum_{j=0}^i (-1)^{n_l-j}\psi^{i-j}(u)\frac{\ddr^{n_l}}{\ddr u^{n_l}}\psi^{j}(u),
    \end{align*}
    where the last identity is obtained by resorting to Hoppe's formula, a convenient variant of the more popular Fa\`a di Bruno's formula \citep[see, e.g.,][]{Joh02}. Using the notation introduced in \eqref{eq:xi} we can write
    \begin{equation}\label{eq:third}
       \E\left[\e^{-u\tilde \mu(\{x_0\})}\tilde \mu^{n_l}(\{x_0\})\right]=\sum_{i=1}^{n_l}c^i\zeta^i \edr^{c\zeta \psi(u)}\xi_{n_l,i}(u).
    \end{equation}
    By plugging \eqref{eq:first}, \eqref{eq:second} and \eqref{eq:third} into \eqref{eq:prod_exp} we get
    \begin{equation}\label{eq:second_case}
        \E\left[ \edr^{-u\tilde\mu(\X)}\prod_{j=1}^k \tilde\mu^{n_j}(\ddr x_j)\right]\cong c^{k-1}(1-\zeta)^{k-1}\edr^{-c\psi(u)}\prod_{m\neq l}\tau_{n_m}(u)P^*(\ddr x_m)\sum_{i=1}^{n_l}c^i\zeta^i\xi_{n_l,i}(u).
    \end{equation}
\end{itemize}
A first-order approximation of the expected value $\E[\exp\{-u\tilde\mu(\X)\}\prod_{j=1}^k\tilde \mu^{n_j}(\ddr x_j)]$ can then be written as linear combination of \eqref{eq:first_case} and \eqref{eq:second_case}, that is
\begin{multline*}
    \left(1-\sum_{l=1}^{k}\delta_{x_0}(\ddr x_l)\right)c^k(1-\zeta)^k\edr^{-c\psi(u)}\prod_{m=1}^k P^*(\ddr x_m)\tau_{n_m}(u)\\
    +c^{k-1}(1-\zeta)^{k-1}\sum_{l=1}^k\delta_{x_0}(\ddr x_l) \left(\prod_{m\neq l}P^*(\ddr x_m)\tau_{n_m}(u)\right)\sum_{i=1}^{n_l}c^i\zeta^i
\xi_{n_l,i}(u).
\end{multline*}
The proof is completed by letting $\varepsilon$ go to 0, by replacing the expected value in \eqref{eq:EPPF_gamma_int} with the last expression, and by computing the integral over $\X^k$. \hfill \qed

\subsection{Proof of Corollary~\ref{thm:interpret_eppf_nrmi}}
 
From Theorem \ref{eppf:nrmi} we obtain the following representation for $\Pi_k^{(n)}(n_1,\ldots,n_k;\zeta)$
\begin{multline*}
\Pi_k^{(n)}(n_1,\ldots,n_k;\zeta) = 
(1-\zeta)^k \Pi_k^{(n)}(n_1,\ldots,n_k;0)+(1-\zeta)^{k-1}\, \sum_{j=1}^k\frac{1}{\Gamma(n-n_j)\Gamma(n_j)}\\[4pt]
\times\: \int_{\X_0^{k-1}}\E \Big[\frac{\tilde\mu^{n_j}(\{x_0\})}{\tilde\mu^{n_j} (\X)}\:\prod_{\ell\ne j}\frac{\tilde\mu^{n_\ell}(\ddr x_\ell)}{\tilde\mu^{n_\ell} (\X)}\Big]
\end{multline*}
where $\X_0^{k-1}=(\X\setminus \{x_0\})^{k-1}$. We do now focus on the integral appearing in the second summand and note that 
\begin{align*}
\int_{\X_0^{k-1}}\E & \Big[\frac{\tilde\mu^{n_j}(\{x_0\})}{\tilde\mu^{n_j} (\X)}\:\prod_{\ell\ne j}\frac{\tilde\mu^{n_\ell}(\ddr x_\ell)}{\tilde\mu^{n_\ell} (\X)}\Big]\\[4pt]
&=\int_0^\infty \int_0^\infty u^{n-n_j-1}\, v^{n_j-1}\: \int_{\X^{k-1}}\E\Big[ \edr^{-(u+v)\,\tilde\mu(\X)} \tilde\mu^{n_j}(\{x_0\})\:\prod_{\ell\ne j}\tilde\mu^{n_\ell}(\ddr x_\ell) \Big]
\:\ddr v\,\ddr u\\[4pt]
&=\int_0^\infty \int_0^\infty u^{n-n_j-1}\, v^{n_j-1}\:\edr^{-c\psi(u+v)}\Big(\prod_{\ell\ne j} \tau_{n_\ell}(u+v)\Big)\\[4pt]
&\qquad \times\: \sum_{i=1}^{n_j} c^i\,\zeta^i\:\frac{1}{i!}\: \sum_{\Delta_{i,n_j}}\binom{n_j}{q_1\,\cdots\, q_i}\:\prod_{r=1}^i \tau_{n_r}(u+v)\,\ddr u\,\ddr v
\end{align*}
where $\Delta_{i,n_j}$ is the set of all vectors of positive integers $(q_1,\ldots,q_i)$ such that $\sum_{r=1}^{i}q_r=n_j$, and we use the fact that 
\[
\xi_{n_j,i}(u+v)=\frac{1}{i!}\sum_{\Delta_{i,n_j}}\binom{n_j}{q_1\cdots q_i}
\prod_{r=1}^i \tau_{q_r}(u+v).
\]
If one recalls the representation of the exchangeable partition probability function in \eqref{eq:eppf:nrmi}, then 
\begin{align*}
&\int_{\X_0^{k-1}}\E  \Big[\frac{\tilde\mu^{n_j}(\{x_0\})}{\tilde\mu^{n_j} (\X)}\:\prod_{\ell\ne j}\frac{\tilde\mu^{n_\ell}(\ddr x_\ell)}{\tilde\mu^{n_\ell} (\X)}\Big]\\[4pt]
& \qquad =\Pi_k^{(n)}(n_1,\ldots,n_{j-1},n_{j+1},\ldots,n_k;0)\  \int_0^\infty \frac{u^{n-n_j-1}\,\edr^{-\psi(u)}\:\prod_{\ell\ne j}\tau_{n_\ell}(u)}{\int_0^\infty u^{n-n_j}\,\edr^{-\psi(u)}\:\prod_{\ell\ne j}\tau_{n_\ell}(u)\:\ddr u}\\[4pt]
&\qquad \qquad  \times \: \int_0^\infty v^{n_j-1}\: \edr^{-c\psi(u+v)}\Big(\prod_{\ell\ne j} \tau_{n_\ell}(u+v)\Big) \sum_{i=1}^{n_j} c^i\,\zeta^i\:\frac{1}{i!}\: \sum_{\Delta_{i,n_j}}\binom{n_j}{q_1\,\cdots\, q_i}\:\prod_{r=1}^i \tau_{n_r}(u+v)\,\ddr u\,\ddr v\\[4pt]
&\qquad =\Pi_k^{(n)}(n_1,\ldots,n_{j-1},n_{j+1},\ldots,n_k;0)\ \int_0^\infty f_j(u)\: \int_0^\infty v^{n_j-1}\: \edr^{-c\psi(u+v)}\Big(\prod_{\ell\ne j} \tau_{n_\ell}(u+v)\Big)\\[4pt]
&\qquad \qquad \times\: \sum_{i=1}^{n_j} c^i\,\zeta^i\:\frac{1}{i!}\: \sum_{\Delta_{i,n_j}}\binom{n_j}{q_1\,\cdots\, q_i}\:\prod_{r=1}^i \tau_{n_r}(u+v)\,\ddr u\,\ddr v
\end{align*}
where, following \cite{Jam09}, 
\begin{equation*}
f_{j}(u)\propto u^{n-n_j-1}\edr^{-c\psi(u)}\prod_{\ell\ne j}\tau_{n_\ell}(u)
\label{eq:density_latent_u}
\end{equation*}
is the density function of a latent random variable $U_{n-n_j}$, conditional on on the $n-n_j$ observations $X_{- j}^{(n)}$ with $k-1$ distinct values that do not include $x_0$ and, hence, are not in cluster $j$. Henceforth we denote these distinct values as $(x^*_{1,j},\ldots,x^*_{k-1,j})$. 
At this point, we benefit from the posterior representation of $\tilde \mu$ given in Theorem 1 of \cite{Jam09}, which entails that, conditional on $X_{- j}^{(n)}$ and on $U_{n-n_j}$, the distribution of $\tilde \mu$ equals the distribution of 
\begin{equation*}
\tilde\mu_j^{(u)}+\sum_{i=1}^{k-1}J_{i,j}^{(u)}\delta_{x^*_{i,j}},
\label{eq:posterior_crm}
\end{equation*}
where $\tilde\mu_j^{(u)}$ is a completely random measure without fixed discontinuities and with intensity $\edr^{-us}\rho(s)\ddr s \,c\, P_0(\ddr x)$ and the jumps $J_{i,j}^{(u)}$ are independent with respective distributions having density $f_{i,j}(s\mid u)\propto s^{n_i}\edr^{-us}\rho(s)$. Hence, a straightforward application of Theorem 2 of \cite{Jam09} shows that, for any $u>0$,
\begin{multline*}
\E\left[\frac{\tilde\mu_j^{(u)}(\{x_0\})}{T_j^{(u)}+\sum_{i}J_{i,j}^{(u)}}\right]^{n_j}=
\int_0^\infty v^{n_j-1}\: \edr^{-c\psi(u+v)}\Big(\prod_{\ell\ne j} \tau_{n_\ell}(u+v)\Big)
\: \sum_{i=1}^{n_j} c^i\,\zeta^i\:\frac{1}{i!}\: \sum_{\Delta_{i,n_j}}\binom{n_j}{q_1\,\cdots\, q_i}\:\prod_{r=1}^i \tau_{n_r}(u+v)\,\ddr v
\end{multline*}
where $T_j^{(u)}=\tilde\mu_j^{(u)}(\X)$. 
Thus we have
\begin{multline}
\label{eq:eppf_ter}
	\Pi_k^{(n)}(n_1,\ldots,n_k;\zeta)=
	(1-\zeta)^k\Pi_k^{(n)}(n_1,\ldots,n_k;0) 
	+(1-\zeta)^{k-1}\sum_{j=1}^k \Pi_{k-1}^{(n-n_j)}(n_1,\ldots,n_{j-1},n_j,\ldots,n_k;0)\\
    \times\: \int_0^\infty f_{j}(u) 
	\E\left[\frac{\tilde\mu_j^{(u)}(\{x_0\})}{T_j^{(u)}+\sum_{i=1}^{k-1}J_{i,j}^{(u)}}\right]^{n_j}\ddr u.
	\end{multline}
The proof is completed upon noting that the integral appearing in right-hand side of \eqref{eq:eppf_ter}
is the conditional probability that $n_j$ observations in $X^{(n)}$ equal $x_0$, given the remaining $n-n_j$ in the sample $X_{-j}^{(n)}$ of $n-n_j$ all differ from $x_0$ and clustered into $k-1$ groups with respective frequencies $n_1,\ldots,n_{j-1},n_{j+1},\ldots,n_k$. 
\hfill \qed

\subsection{Proof of Theorem \ref{thm:K_n}}

The proof follows by combining \eqref{eq:eppf_bis} in Theorem \ref{thm:interpret_eppf_nrmi} with the fact that 
\begin{equation*}
    \Pr(K_n=k;\zeta)=\frac{1}{k!}\sum_{\Delta_{k,n}}\binom{n}{q_1 \,\cdots q_k} \Pi_k^{(n)}(q_1,\ldots,q_k;\zeta).
\end{equation*}
\hfill \qed

\subsection{Proof of Theorem \ref{thm:predictive}}

For the purpose of the proof, we introduce the quantity 
\begin{equation*}
    \Pi_{k,j}^{(n)}(n_1,\ldots,n_k;\zeta),
\end{equation*}
for $j=0,1,\ldots,k$,
to denote the probability of observing a partition of $n$ observations into $k$ distinct blocks such that the frequency of the block coinciding with $x_0$ is given by the $j$-th argument of the function. The case $j=0$ indicates that no block in the partition coincides with $x_0$. From \eqref{eq:eppf:nrmi} we have that 
\begin{align*}
    \Pi_{k,0}^{(n)}(n_1,\ldots,n_k;\zeta)&=\frac{1}{\Gamma(n)}c^k (1-\zeta)^k \int_0^\infty u^{n-1}\e^{-c\psi(u)}\prod_{m=1}^k \tau_{n_m}(u)\d u,\\
    \Pi_{k,j}^{(n)}(n_1,\ldots,n_k;\zeta)&=\frac{1}{\Gamma(n)}c^{k-1} (1-\zeta)^{k-1} \sum_{i=1}^{n_j}c^i\zeta^i \int_0^\infty u^{n-1}\e^{-c\psi(u)} \xi_{n_j,i}(u)\prod_{m\neq j}^k \tau_{n_m}(u)\d u.
\end{align*}
The predictive distribution takes the form 
\begin{equation}\label{eq:pred_general}
    \Pr(X_{n+1}\in A\mid X^{(n)};\zeta)=w_{k,n}^{(*)} P^*(A)+w_{k,n}^{(0)}\delta_{x_0}(A)+\sum_{l:x^*_l\neq x_0} w_{k,n}^{(l)} \delta_{x_l^*}(A),
\end{equation}
where the weights in \eqref{eq:pred_general} are as follows. If $x_0\notin\{x_1^*,\ldots,x_k^*\}$,
\begin{align*}
    w_{k,n}^{(*)}&=\frac{\Pi_{k+1,0}^{(n+1)}(n_1,\ldots,n_k,1;\zeta)}{\Pi_{k,0}^{(n)}(n_1,\ldots,n_k;\zeta)},\\
        w_{k,n}^{(0)}&=\frac{\Pi_{k+1,k+1}^{(n+1)}(n_1,\ldots,n_k,1;\zeta)}{\Pi_{k,0}^{(n)}(n_1,\ldots,n_k;\zeta)},\\
        w_{k,n}^{(l)}&=\frac{\Pi_{k,0}^{(n+1)}(n_1,\ldots,n_l+1,\ldots,n_k;\zeta)}{\Pi_{k,0}^{(n)}(n_1,\ldots,n_l,\ldots,n_k;\zeta)}.
\end{align*}
On the other hand, if $x_0=x_j^*$ for some $j=1,\ldots,k$, then
\begin{align*}
    w_{k,n}^{(*)}&=\frac{\Pi_{k+1,j}^{(n+1)}(n_1,\ldots,n_k,1;\zeta)}{\Pi_{k,j}^{(n)}(n_1,\ldots,n_k;\zeta)},\\
        w_{k,n}^{(0)}&=\frac{\Pi_{k,j}^{(n+1)}(n_1,\ldots,n_j+1,\ldots,n_k;\zeta)}{\Pi_{k,j}^{(n)}(n_1,\ldots,n_j,\ldots,n_k;\zeta)},\\
        w_{k,n}^{(l)}&=\frac{\Pi_{k,j}^{(n+1)}(n_1,\ldots,n_l+1,\ldots,n_k;\zeta)}{\Pi_{k,j}^{(n)}(n_1,\ldots,n_l,\ldots,n_k;\zeta)}.
\end{align*}
Simple algebra completes the proof.
\hfill \qed

\subsection{Proof of Theorem \ref{thm:N0_norm_ss}}
As in the proof of Theorem \ref{thm:predictive}, we use the notation
$$\Pi_{k,j}^{(n)}(n_1,\ldots,n_k;\zeta),\quad j\in\{0,1,\ldots,k\}$$
to indicate the probability that an inner spike and slab hNRMI model induce a partition of $n$ observations into $k$ blocks of size $n_1,\ldots,n_k$, where, if $j\in\{1,\ldots,k\}$, the observations of the $j$-th block coincide with $x_0$, while if $j=0$ then no observation coincides with $x_0$. Along similar lines as Theorem \ref{eq:eppf:nrmi}, we can show that 
\begin{itemize}
    \item if $j=0$
    \begin{align*}
    \Pi_{k,0}^{(n)}(n_1,\ldots,n_k;\zeta)=\frac{1}{\Gamma(n)}c^k(1-\zeta)^k \int_0^\infty u^{n-1}\edr^{-c\psi(u)}\prod_{m=1}^k \tau_{n_m}(u)\ddr u;
    \end{align*}
    \item if $j\in\{1,\ldots,k\}$
    \begin{align*}
        \Pi_{k,j}^{(n)}(n_1,\ldots,n_k;\zeta)=\frac{1}{\Gamma(n)}c^{k-1}(1-\zeta)^{k-1} \sum_{i=1}^{n_j}c^i\zeta^i \int_0^\infty u^{n-1}\edr^{-c\psi(u)}\xi_{n_j,i}(u)\prod_{m\neq j}\tau_{n_m}(u)\ddr u.
    \end{align*}
\end{itemize}
Next we observe that 
\begin{align}
    \Pr(N_0^{(n)}=0)&=\sum_{k=1}^n \frac{1}{k!}\sum_{\bm{n}\in \mathcal{N}_{k}^{(n)}}\binom{n}{n_1 \cdots n_k}\Pi_{k,0}^{(n)}(n_1,\ldots,n_k;\zeta)\notag\\
    &=\frac{1}{\Gamma(n)}\sum_{k=1}^n c^k (1-\zeta)^k \int_0^\infty u^{n-1}\edr^{-c\psi(u)}\frac{1}{k!}\sum_{\bm{n}\in\mathcal{N}_{k}^{(n)}}\binom{n}{n_1 \cdots n_k}\prod_{m=1}^k \tau_{n_m}(u)\ddr u.\label{eq:Pr_N0_0}
\end{align}
Similarly, if $j\in\{1,\ldots,n\}$,
\begin{align}
    \Pr(N_0^{(n)}=j)&=\sum_{k=1}^{n-j+1}\frac{1}{(k-1)!} \sum_{\bm{n}\in\mathcal{N}_{k-1}^{(n-j)}} \binom{n}{n_1 \cdots n_{k-1}}k\, \Pi_{k,k}^{(n)}(n_1,\ldots,n_{k-1},j;\zeta)\notag\\
    &=\frac{1}{\Gamma(n)} \sum_{k=1}^{n-j+1} k\, c^{k-1}(1-\zeta)^{k-1}\notag \\
    &\times \sum_{i=1}^j c^i \zeta^i \int_0^{\infty} u^{n-1}\edr^{-c\psi(u)}\xi_{j,i}(u)\frac{1}{(k-1)!}\sum_{\bm{n}\in\mathcal{N}_{k-1}^{(n-j)}}\binom{n}{n_1\cdots n_{k-1}}\prod_{m=1}^{k-1}\tau_{n_m}(u)\ddr u\label{eq:Pr_N0_j},
    \end{align}
with $\sum_{\bm{n}\in\mathcal{N}_{k-1}^{(n-j)}} \equiv 1$ when $k=1$. The two expressions in \eqref{eq:Pr_N0_0} and \eqref{eq:Pr_N0_j} can be summarized in one, leading to \eqref{eq:N0_norm_ss},
where we observe that $\xi_{j,0}=0$ if $j\geq 1$, and $\xi_{0,0}=1$.
\hfill \qed 

\bibliographystyle{elsarticle-harv}
\bibliography{refs.bib}

\end{document}